\documentclass[final,5p,times,twocolumn]{elsarticle}

\usepackage{amssymb}
\usepackage{amsmath}
\usepackage{bbm}
\usepackage{graphicx}
\usepackage{exscale}
\usepackage{clrscode}
\usepackage{xcolor}
\usepackage{color, colortbl}
\usepackage{booktabs}

\biboptions{sort&compress}

\newcommand{\E}[1]{\mathop{{\rm \bf E}\!\left\{#1\right\}}\nolimits}

\newdefinition{dfn}{Definition}

\begin{document}

\begin{frontmatter}
\title{Estimation of market efficiency process within time-varying autoregressive models by extended Kalman filtering approach}

\author[CEMAT]{M.~V.~Kulikova\corref{cor}} \ead{maria.kulikova@ist.utl.pt} \cortext[cor]{Corresponding
author.}

\author[CEMAT]{G.~Yu.~Kulikov} \ead{gkulikov@math.ist.utl.pt}

\address[CEMAT]{CEMAT, Instituto Superior T\'ecnico, Universidade de Lisboa, Av.~Rovisco Pais, 1049-001 Lisboa, Portugal.}

\begin{abstract}
This paper explores a time-varying version of weak-form market efficiency
that is a key component of the so-called Adaptive Market Hypothesis (AMH). One of the most common methodologies used for modeling and estimating a degree of market efficiency lies in an analysis of the serial autocorrelation in observed return series. Under the AMH, a time-varying market efficiency level is modeled by time-varying autoregressive (AR) process and traditionally estimated by the Kalman filter (KF). Being a linear estimator, the KF is hardly capable to track the hidden nonlinear dynamics that is an essential feature of the models under investigation. The contribution of this paper  is threefold. We first provide a brief overview of time-varying AR models and estimation methods utilized for testing a weak-form market efficiency in econometrics literature. Secondly, we propose novel accurate estimation approach for recovering the hidden process of evolving market efficiency level by the extended Kalman filter (EKF). Thirdly, our empirical study concerns an examination of the Standard and Poor's 500 Composite stock index and the Dow Jones Industrial Average index.
Monthly data covers the period from November 1927 to June 2020, which includes the U.S. Great Depression, the 2008-2009 global financial crisis and the first wave of recent COVID-19 recession. The results reveal that the U.S. market was affected during all these periods, but generally remained weak-form
efficient since the mid of 1946 as detected by the estimator.
\end{abstract}


\end{frontmatter}

\section{Introduction}\label{sect1}

The Adaptive Market Hypothesis (AMH) concept with the assumption of possibly time-varying weak-form market efficiency has received an increasing attention in recent years. Its popularity has spread rapidly since the mid-2000s and the first papers promoted this idea in the financial industry~\cite{2004:Lo,2005:Lo}. An excellent and systematic review of weak-form market efficiency literature can be found in~\cite{2011:Lim}. Apart from empirical studies, a large number of research papers has been aimed to model and estimate the dynamics of time-varying market efficiency level from the return history available~\cite{1997:Emerson,1999:Zalewska,2000:Rockinger,2003:Li:AMH:China,2008:Posta,2014:Ito,2016:Ito,2016:Noda,2019:RJ:Kulikov}.

One of the most common methodologies used for modeling and estimating a degree of market efficiency lies in an analysis of the serial autocorrelation in observed return series. It follows from the definition of weak-form market efficiency introduced in~\cite{1965:Fama}: ``a market is weak-form efficient when there is no predictable profit opportunity based on the past movement of asset prices''. This implies an evident conclusion: if the existence of a serial autocorrelation in return series is observed, then some degree of market inefficiency appears. Consequently, a common approach to modeling and tracking the (time-varying) level of market efficiency within the AMH framework is to apply an autoregressive (AR) model with time-dependent coefficients
to a chosen history of returns~\cite{2009:ItoSugiyama,2011:Kim:AMH}. We may mention that time-varying autoregressive models are often utilized in econometric literature along with Bayesian estimation methods. For instance, the readers are referred to a recent research devoted to Bayesian compressed vector autoregression for financial time-series analysis and forecasting in~\cite{2019:Taveeapiradeecharoen,2020:Aunsri} and many other works. In this paper, we focus on the tests for {\it evolving market efficiency} based on time-varying AR models, only.

A simple approach to design the tests under the AMH is to assume that the return series follows a time-varying AR model with {\it homoscedastic conditional variance} assumption, see, for example,~\cite{2009:ItoSugiyama}. However, modern econometric trends suggest the use of more sophisticated system structures to model complicated `stylized' facts such as volatility clustering via Autoregressive Conditional Heteroscedasticity (ARCH) and generalized ARCH (GARCH) models~\cite{1982:Engle,1986:Bollerslev}. The tests for a weak-form  market efficiency based on time-varying AR models
with GARCH process for modeling dynamics of conditional variance process are still valid because an uncorrelated time series might exhibit autocorrelation in the squared series. We also refer to the discussion in econometric literature of a difference between the `pure white noise series' (i.e. a random walk with independent and identically distributed  (i.i.d.)  increments) are `white noise' (i.e. serially uncorrelated) and their appropriate utilization in the statistical tests for time-varying market efficiency in~\cite[Section~3.3]{2011:Lim} as well as statistical  tools developed for uncovering  hidden  nonlinear structures in previously observed serially uncorrelated stock market in~\cite{1985:Hinich,1988:Brockett,1989:DeGooijer,1989:Scheinkman}.

Tests for a weak-form  market efficiency based on time-varying AR models with various GARCH-type processes accounting for {\it heteroscedasticity} in the conditional variance have received increased attention in recent years. For instance, the test for evolving efficiency (TEE) proposed in~\cite{1997:Emerson} and formalized in~\cite{1999:Zalewska} employs a GARCH-in-Mean(1,1) process combined with an AR(1) assumption to model the conditional mean where the regression coefficients follows a random walk and should be estimated when formulating a conclusion about the market's efficiency. Empirical studies based on the TEE approach include seven African stock markets in~\cite{2005:Jefferis}, eleven Arab stock markets in~\cite{2010:Abdmoulah}, as well as the Gulf Cooperation Council (the economic union consisting of all Arab states of the Persian Gulf, except Iraq) in~\cite{2016:Charfeddine}, the Russian and Czech markets in~\cite{2002:Hall,2008:Posta,2019:RJ:Kulikov}, and so on. In~\cite{2000:Rockinger,2003:Li:AMH:China}, asymmetric and threshold GARCH-type specifications  are proposed in evolving market efficiency tests that allow for asymmetries in the volatility of returns in reaction to information shocks. This is referred to as the leverage effect, and is often modeled using asymmetric GARCH-type processes. For example, the empirical study in~\cite{2003:Li:AMH:China} found significant evidence of a leverage effect in the Shanghai stock market returns while the Shenzhen return series show no asymmetry.

Overall, the main drawback of the previously suggested tools for market efficiency research under the AMH is an application of the classical Kalman filter (KF) to calibrate the models and to extract the hidden evolving efficiency process. Being a linear estimator, the KF is hardly capable to track the hidden nonlinear dynamics that is an essential feature of the models under investigation.
The contribution of this paper  is threefold. We first provide a brief overview of time-varying AR models and estimation methods utilized for testing a weak-form market efficiency in econometrics literature. Secondly, we propose novel accurate estimation approach for recovering the hidden process of evolving market efficiency level by the extended Kalman filter (EKF). The proposed methodology solves two estimation problems simultaneously that is the model calibration and the unknown state estimation that allows for tracking and detecting the inefficiency periods (i.e. anomalies) from return series.  Thirdly, our empirical study concerns an examination of the Standard and Poor's 500 Composite stock index and the Dow Jones Industrial Average index.
Monthly data covers the period from November 1927 to June 2020, which includes the U.S. Great Depression, the 2008-2009 global financial crisis and the first wave of recent COVID-19 recession. The results reveal that the U.S. market was affected during all these periods, but generally remained weak-form
efficient since the mid of 1946 as detected by the estimator.

The rest of the paper is organized as follows. Major modeling approaches under the state-space representation and related filtering methods utilized for recovering the evolving market efficiency level are summarized in Section~\ref{survey:models}. Moving window tests are applied for time-varying
sample autocorrelation estimation and for the analysis of changing serial dependence over time for U.S. stock market in
Section~\ref{sec:data}. The nonlinear estimation methods for recovering the dynamic market efficiency level for models with
both homoskedastic and heteroscedastic conditional variance assumptions are proposed in Section~\ref{Section:models}. Finally, Section~\ref{Section:conclusion}  summarizes the key findings of our study
and concludes the paper. The details of the extended Kalman filtering approach and the parameters estimation procedure within the method of maximum likelihood are presented in Appendix.

\section{Survey of the AMH state-space modeling techniques} \label{survey:models}

A common approach to modeling and tracking the level of market efficiency
is to apply an autoregressive (AR) model to a chosen history of returns for estimating its coefficients.
The AR process of order $n$ is given as follows:
\begin{equation}
y_t = \sum \limits_{i=1}^{n} \beta_i y_{t-i} + \varepsilon_t, \quad \varepsilon_t \sim IID(0,\sigma_{\varepsilon}^2)  \label{AR:n}
\end{equation}
where $y_t$ is the log return series, $\beta_i$, $i=1,\ldots,n$ are the parameters of the model, which are assumed to be constant over time within the classical EMH concept.

Under the EMH, tests for weak-form  market efficiency based on an AR($n$) modeling  approach yield the following condition: $\beta_1=\beta_2= \ldots =\beta_n=0$. If this condition is met, equation~\eqref{AR:n} takes the form of $y_t = \varepsilon_t$, which implies that the returns are independently distributed and that, consequently, the historical price information cannot provide consistent profit opportunities. This corresponds to the definition of weak-form market efficiency in~\cite{1965:Fama}. In fact, the assumption $\varepsilon_t \sim IID(0,\sigma_{\varepsilon}^2)$ is too strong for testing weak-form market efficiency. Following the discussion in~\cite{2011:Lim}, it is sufficient to assume that $\varepsilon_t$ is a white noise process. If $\beta_1=\beta_2= \ldots =\beta_n=0$, then this implies that the returns also follow a white noise process. This means that they are serially uncorrelated in time, and
again that the historical price information cannot provide profit opportunities, i.e. the market is
weak-form efficient by definition. Finally, model~\eqref{AR:n} may include a constant trend $\beta_0$, but the test for market efficiency remains the same $\beta_1=\beta_2= \ldots =\beta_n=0$. Indeed, this condition yields $y_t = \beta_0 + \varepsilon_t$ where $\varepsilon_t$ is a white noise process, i.e. the log return is a white noise process about the mean and, hence, the returns are still serially uncorrelated in time.

Within the AMH approach, a common technique to modeling and tracking the {\it time-varying} level of market efficiency
is to apply an autoregressive (AR) model with {\it time-dependent} coefficients to a return history. Here we explore
only AR models in the presence of Gaussian uncertainties as suggested in~\cite{1997:Dahlhaus}:
\begin{align}
y_t & = \sum \limits_{i=1}^{n} \beta_{i,t} y_{t-i} + \varepsilon_t, & \varepsilon_t & \sim {\cal N}(0,\sigma_{\varepsilon}^2) \label{AR:n:time}
\end{align}
where $\beta_{i,t}:=\beta_i(t_k)$ and $\varepsilon_t$ is a Gaussian white noise process with zero mean and variance $\sigma_{\varepsilon}^2>0$. The market efficiency dynamics is reflected in the time-varying regression coefficients $\beta_{i,t}$, $i=1, \ldots, n$, which should be modeled and estimated.

A random walk specification is the most commonly utilized model for time-varying regression coefficients employed in the econometric literature dedicated to evolving market efficiency research. The readers are also referred to a brief survey of modeling techniques presented in Table~\ref{tab:models}. A random walk specification for time-varying regression coefficients yields
\begin{align}
\beta_{i,t} & = \beta_{i,t-1} + w_{i,t}, \; i=1,\ldots,n, & w_{i,t} & \sim {\cal N}(0,\sigma_{w_i}^2) \label{Beta:RW}
\end{align}
where $\sigma_{w_i}^2 \ge 0$, $i=1,\ldots,n$. The disturbances $\varepsilon_t$ and $w_{i,t}$ are mutually uncorrelated white noise processes.

\begin{table*}[ht!]
{\scriptsize
\caption{A brief overview of time-varying autoregressive models and estimation methods utilized for testing a weak-form market efficiency.} \label{tab:models}                                                                                                                                                                                                                                                                                                                                                                                                                                                                                                                                                                                                                                                                \begin{tabular}{llcccccl}
\toprule %\rowcolor{black!5}
{\bf Year} & {\bf Who and} & {\bf Model}     & {\bf Model}  & {\bf State and} & {\bf Parameters and} & {\bf Efficiency} & {\bf Applied to} \\
           & {\bf Where}   & {\bf Equations} & {\bf Description} & {\bf Estimator} &  {\bf Estimator}     & {\bf Criterion}  &  {\bf Index and Period}\\
\toprule
 1997 & Emerson et al.  & $y_t = \beta_{0,t} + \sum \limits_{i=1}^{p=2} \beta_{i,t} y_{t-i} + \delta h_t + \varepsilon_t, $ & GARCH-M(1,1) with
 & $[h_t, \beta_{0,t}, \beta_{1,t}, \beta_{2,t}]$  & $\delta, a_0, a_1, b_1, \sigma_{i,w}^2$ & $\hat \beta_{i,t} = 0$,    & Sofia Stock Exchange \\
 &  in~\cite{1997:Emerson} &  $\varepsilon_t \sim {\cal N}(0,h_t), h_t =  a_0 + a_1 \varepsilon_{t-1}^2 + b_1 h_{t-1}$ & AR(2) for the mean eq.& by KF &
 $i=\overline{0,2}$ by MLE & $i>0$  &  (Jan.1994 -- Jan.1996)   \\
     &  &  $\beta_{i,t} = \beta_{i,t-1} + w_t, w_t \sim {\cal N}(0,\sigma_{i,w}^2)$ & coeff. follow random walk &  &  &  &\\[10pt]
%______
 1999 & Zalewska-Mitura  & $y_t = \beta_{0,t} + \beta_{1,t} y_{t-1} + \delta h_t + \varepsilon_t$ & GARCH-M(1,1) with
 & $[h_t, \beta_{0,t}, \beta_{1,t}]$  & $\delta, a_0, a_1, b_1, \sigma_{0,w}^2$ & $\hat \beta_{1,t} = 0$    &  Budapest Stock Exchange \\
 &  et al. in~\cite{1999:Zalewska} &  $\varepsilon_t \sim {\cal N}(0,h_t), h_t =  a_0 + a_1 \varepsilon_{t-1}^2 + b_1 h_{t-1}$ & AR(1) for the mean eq.& by KF &
 and $\sigma_{1,w}^2$ by MLE & & {FTSE 100} (U.K.)  \\
     &  &  $\beta_{i,t} = \beta_{i,t-1} + w_t, w_t \sim {\cal N}(0,\sigma_{i,w}^2)$ & coeff. follow random walk &  &  &  & (Jan.1991 -- Oct.1997) \\[10pt]
%______
 2000 & Rockinger  & $y_t = \beta_{0,t} + \beta_{1,t} y_{t-1} + \varepsilon_t, \varepsilon_t \sim {\cal N}(0,h_t)$ &  asymmetric GARCH
 & $[\beta_{0,t}, \beta_{1,t}]$  & $a_0, a_1, a_2, b_1, \sigma_{0,w}^2$ & $\hat \beta_{1,t} = 0$    &  Budapest Stock Exchange \\
 &  et al. in~\cite{2000:Rockinger} &  $h_t =  a_0 + a_1 \varepsilon_{t-1}^2\Lambda_{y_{t-1}>0} $ & AR(1) for the mean eq.& by KF &
 and $\sigma_{1,w}^2$ by MLE & & Czech, Russian indices   \\
     &  &  $+ a_2 \varepsilon_{t-1}^2\Lambda_{y_{t-1}<0} + b_1 h_{t-1}$ & coeff. follow random walk &  &  &  & Warsaw Stock Exchange \\
     &  &  $\beta_{i,t} = \beta_{i,t-1} + w_t, w_t \sim {\cal N}(0,\sigma_{i,w}^2)$ &  &  &  &  & (Apr.1994 -- Jun.1999) \\[10pt]
 %______
 2003 &  Xiao-Ming Li  & $y_t = \beta_{0,t} + \sum \limits_{i=1}^{p=3} \beta_{i,t} y_{t-i} + \varepsilon_t, \varepsilon_t \sim {\cal N}(0,h_t)$ & TGARCH(1,1) with
 & $[\beta_{0,t}, \beta_{1,t}, \beta_{2,t}, \beta_{3,t}]$  & $a_0, a_1, a_1^+, b_1, \sigma_{i,w}^2$ & $\hat \beta_{i,t} = 0$,    & Shenzhen Stock Exchange \\
 &  in~\cite{2003:Li:AMH:China} &  $h_t =  a_0 + a_1 \varepsilon_{t-1}^2 +a_1^+ (\varepsilon^+_{t-1})^2+b_1 h_{t-1}$ & AR(3) for the mean eq.& by KF &
 $i=\overline{0,3}$ by MLE & $i>0$  & (Apr.1991 -- Jan.2001)  \\
  &  &  with $\varepsilon^+_{t-1} = \max\{\varepsilon_{t-1},0\}$  & coeff. follow random walk &  &  &  & Shanghai Stock Exchange \\
     &  &  $\beta_{i,t} = \beta_{i,t-1} + w_t, w_t \sim {\cal N}(0,\sigma_{i,w}^2)$ &  &  &  &  &  (Jan.1991 -- Jan.2001)\\[10pt]
%______
 2008 & Po{\v{s}}ta  & $y_t = \beta_{0} + \beta_{1,t} y_{t-1} +\varepsilon_t$ or & GARCH(1,1) and (2,2),
 & $[h_t, \beta_{1,t}]$  & $\beta_0, a_0, a_1, b_1, \eta$ & $\hat \beta_{1,t} = 0$    &  Prague Stock Exchange \\
 &  in~\cite{2008:Posta} & $y_t = \beta_{0} + \beta_{1,t} y_{t-1} + \delta h_t + \varepsilon_t$   &  GARCH-M(1,1),  & by KF &
 etc. by MLE & & (Jan.1995 -- Jul.2007)  \\
     &  &  $\varepsilon_t \sim {\cal N}(0,h_t)$ where & E-GARCH, TARCH, &  &  &  &  \\
     &  &  $h_t$ follows various GARCH & AR(1) for the mean eq.  &   &  &  &  \\
     &  &  $\beta_{1,t} = \beta_{1,t-1} + w_t, w_t \sim {\cal N}(0,e^{\eta})$ & coeff. follow random walk &  &  &  &  \\[10pt]
%______
2009 & Ito et al. & $y_t = \beta_{1,t} y_{t-1} + \varepsilon_t, \varepsilon_t \sim {\cal N}(0,\sigma_{\varepsilon}^2)$ & AR(1) with coefficients     & $\beta_{1,t}$ by KF & $\sigma_{\varepsilon}^2$, $\sigma_{w}^2$ & $\hat \beta_{1,t} = 0$   & S\&P500 (U.S.) \\
     & in~\cite{2009:ItoSugiyama} &  $\beta_{1,t} = \beta_{1,t-1} + w_t, w_t \sim {\cal N}(0,\sigma_{w}^2)$ & follow random walk & smoother in~\cite{2007:Ito:note}  &  &  & (Jan.1955 -- Feb.2006) \\[10pt]
%___
 2019 & Kulikov et al.  & $y_t = \beta_{0,t} + \beta_{1,t} y_{t-1} + f(\delta) h_t + \varepsilon_t, $ & GARCH-M(1,1) with
 & $[h_t, \beta_{0,t}, \beta_{1,t}]$  & $\delta, a_0, a_1, b_1, \sigma_{0,w}^2$ & $\hat \beta_{1,t} = 0$    & {FTSE 100} (U.K.) \\
 &  in~\cite{2019:RJ:Kulikov} &  $h_t =  a_0 + a_1 \varepsilon_{t-1}^2 + b_1 h_{t-1}$ or & stochastic GARCH-M  & by EKF &
 $\sigma_{1,w}^2$, $\sigma_{u}^2$ by MLE & & {TOP 40} (S.A.)   \\
     &  &  $h_t =  a_0 + a_1 \varepsilon_{t-1}^2 + b_1 h_{t-1}+u_t$ & AR(1) for the mean eq. &  &  &  & {NSE 20} (Kenya) \\
     &  &  $\varepsilon_t \sim {\cal N}(0,h_t)$, $u_t \sim {\cal N}(0,\sigma_{u}^2)$ &  and nonlinear feedback &  &  &  & {RTSI} (Russia) \\
     &  &  $\beta_{i,t} = \beta_{i,t-1} + w_t, w_t \sim {\cal N}(0,\sigma_{i,w}^2)$ & coeff. follow random walk &  &  &  & (Mar.2002 -- Mar.2006) \\
\bottomrule
\end{tabular}}
\end{table*}

In summary, the most simple test for {\it evolving} weak-form market efficiency consists of equations~\eqref{AR:n:time}, \eqref{Beta:RW} and implies the time-dependent AR coefficient estimation procedure. Following Table~\ref{tab:models}, this model specification has been used for empirical study of the S\&P500 index and its efficiency level estimation in~\cite{2009:ItoSugiyama}. The same model has been later used for tracking the evolution of return predictability of the U.S. stock market based on the Dow Jones Industrial Average index in~\cite{2011:Kim:AMH}. Under the AMH, the tests imply a detection procedure of the market efficiency periods when $\beta_{1,t}=\beta_{2,t}= \ldots =\beta_{n,t}=0$. We need to emphasize the fact that first-order autoregressive models, AR(1), are traditionally used in practice because they are rich enough to be interesting and simple enough to
permit complicated extensions for modeling the coefficients' dynamics.

One of the most important and sophisticated extensions mentioned above is to model the {\it heteroscedasticity}
frequently observed in return series and integrated in the tests for evolving weak-form efficiency. More precisely, in contrast to a constant variance $\sigma_{\varepsilon}^2$ assumption in equation~\eqref{AR:n:time}, a more realistic modeling approach suggests to take into account the time-varying volatility process observed in return series. The most demanded modeling strategy, which takes into account this fact, yields {\it the test for evolving efficiency} (TEE) proposed at the first time in~\cite{1997:Emerson} and rationalized afterwards in~\cite{1999:Zalewska}. A popularity of the TEE methodology used in empirical study of various stock markets has been spread rapidly since the beginning of the 2000s. The TEE implies a GARCH-in-Mean(1,1) process combined with AR(1) equation for modeling the conditional mean. Meanwhile, the TEE regression coefficients are assumed to be time-varying and should be estimated for making a decision about the market efficiency regime. Following Table~\ref{tab:models}, the TEE regression coefficients are still assumed to follow the random walk as well as in all other studies. Recently, the TEE model has been generalized for assessing the impact of a presence of nonlinear volatility feedback in~\cite{2019:RJ:Kulikov}.

Table~\ref{tab:models} briefly summarizes the key time-varying market efficiency modeling strategies within the state-space framework along with filtering methods utilized for hidden state and system parameters estimation. In other words, the filters are used for tracking the hidden process that is a time-varying level of market efficiency under the AMH. As can be seen, the next important step in modeling the time-varying nature of weak-form market efficiency within the AMH methodology has been taken in~\cite{2000:Rockinger}. There, the asymmetric effects on volatility have been taken into account by developing the first test for evolving market efficiency within the {\it threshold} heteroskedastic models. The threshold GARCH  (TGARCH) specifications account for asymmetries in the stock price volatility reaction to information shocks.  The model proposed in~\cite{1992:Campbell,1993:Glosten,1994:Zakoian} assumes $\Lambda_{y_{t-1}>0}$ term that represents a dummy variable, taking the value $1$ when the last period's error $y_{t-1}$ is positive, and taking the value $0$ otherwise. Similarly $\Lambda_{y_{t-1}<0}$
takes the value $1$ if the last period's error is negative, and $0$ otherwise. In summary, the related model for tracking the market efficiency changes over time consists of the TGARCH specification used for conditional variance equation along with AR(1) process utilized for modeling  the conditional mean. Again, the time-varying regression coefficients are assumed to follow a random walk.

An alternative asymmetric GARCH modeling approach, which is capable to model the asymmetric (leverage) effect observed in real return series, has been explored in~\cite{2003:Li:AMH:China}. There, the proposed
model simultaneously accounts for the leverage effect, examines the possibility of information
transmission, and tests for evolving market efficiency. For this, the conditional variance of returns
is assumed to evolve according to an alternative asymmetric GARCH (A-GARCH) process with
an additional term  $(\varepsilon^+_{t-1})^2$  where $\varepsilon^+_{t-1} = \max\{\varepsilon_{t-1},0\}$. This conditional variance equation is further combined with the time-varying regression models. Significant asymmetry, i.e. the present of the leverage effect, can be also interpreted as a consequence of a presence of government's intervention in the investors' activity. The empirical study provided in~\cite{2003:Li:AMH:China} has revealed a significant sign of the leverage effect for the Shanghai stock returns while the Shenzhen return series shows no such asymmetry.

Finally, the asymmetric GARCH models have been investigated in~\cite{2008:Posta}, as well. For instance, it has been found that the E-GARCH specification used for the conditional variance process and combined with the time-varying AR(1) model for the conditional mean is the most appropriate way for describing the market efficiency process of the Prague Stock Exchange.

Having analyzed the filtering methods utilized for estimating the unknown models'parameters and hidden time-varying efficiency process in Table~\ref{tab:models}, we conclude that the classical Kalman filter (KF) is traditionally applied to estimate the state vector. Meanwhile the method of maximum likelihood (MLE) is used for the chosen model calibration. The major drawback of this approach is that the underlying nonlinear dynamics is hardly recovered by the linear estimator, which is the classical KF, in an adequate way.

\section{Data and empirical evidences from U.S. stock market}\label{sec:data}

We explore the U.S. stock market by examining the US Standard and Poor's Composite stock index (S\&P500) along with the Dow Jones Industrial Average index (DJIA). The data is collected on a monthly frequency basis for a period from November 1927 to July 2020. The monthly log-returns are computed in the usual way on a continuously compounded basis, i.e. by taking the log first difference of the time series at hand $y(t_k) = \ln S(t_k) - \ln S(t_{k-1})$ where $S(t_k)$ is the closing price at time instance $t_k$. We note that the examination of daily log-returns often
yields a conclusion of weak-form inefficiency because the tests examine serial autocorrelation,
which is prevalent in daily returns. For this reason, monthly observations are preferable for an
appropriate study of evolving market efficiency; e.g., see the discussion in~\cite{2009:ItoSugiyama}.

\begin{table}[ht!]
{\small
\caption{Summary statistics of the S\&P500 and DJIA monthly log returns} \label{tab:stats}                                                                                                                                                                                                                                                                                                                                                                                                                                                                                                                                                                                                                                                                \begin{tabular}{lrr}
\toprule
Period                 &  \multicolumn{2}{r}{Nov. 1927 -- June 2020}  \\
Return Series &  S\&P500 & DJIA \\
Number of observations & $N=1112$ & $N=1112$\\
\toprule
Mean, $\hat \mu$   &   0.00476 & 0.00440\\
Median &   0.00907 & 0.00857 \\
Standard Deviation, $\hat \sigma$  &  0.05402 &  0.05276\\
Skewness &  -0.62185 &  -0.81894\\
Excess Kurtosis &  7.32653 & 7.02294\\
\hline
$\hat \rho_1$  & 0.0832 & 0.0730 \\
$\hat \rho_{10}$  & 0.0225 &   -0.0019 \\
$\hat \rho_{15}$  & 0.0034 &  0.0219 \\
\hline
$Q(1)$ &  7.7268  &  5.9408 \\
$p$-value & 0.0054 &  0.0148\\[5pt]
$Q(10)$    &  31.9632 &  29.9261 \\
$p$-value &  0.0004 &  0.0009 \\[5pt]
$Q(15)$ &  40.5103  &  35.8616 \\
$p$-value &  0.0004 &  0.0019 \\
\hline
Ljung-Box Q-tests & $H_0$ is rejected & $H_0$ is rejected  \\
up to lag $l=15$      & for some $l$      & for some $l$  \\
\bottomrule
\end{tabular}}
\end{table}

The summary statistics for the monthly return series under examination can be found in Table~\ref{tab:stats}. There, the term $\hat \rho_l$ stands for the sample autocorrelation coefficient at lag $l$. Meanwhile, the value $Q(l)$ is the corresponding Ljung-Box statistics in the corresponding Q-test for residual autocorrelation.  This is a joint test for the hypothesis that the first $l$ autocorrelation coefficients are equal to zero. We have tested the samples up to lag $15$ and conclude that the null hypothesis that the residuals are uncorrelated is rejected at both the 1\% and 5\% significance levels for both the S\&P500 and DJIA series. This means that there exists evidence of serial autocorrelation in the monthly returns. Consequently, we have an initial indication that the U.S. market is weak-form inefficient in the {\it absolute sense} over the entire time period examined. To explore {\it evolving} efficiency, we next provide the empirical study based on a rolling window test.

\begin{figure*}[th!]
\begin{tabular}{cc}
\includegraphics[width = 0.5\textwidth]{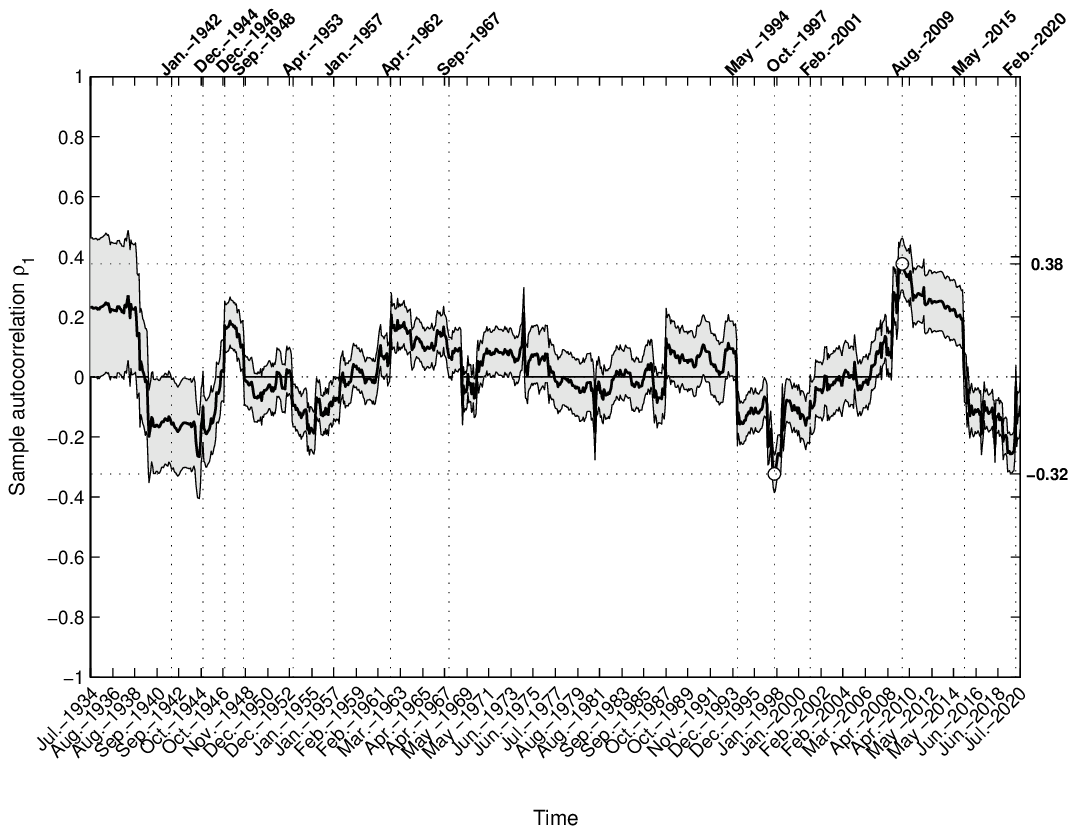} &
\includegraphics[width = 0.5\textwidth]{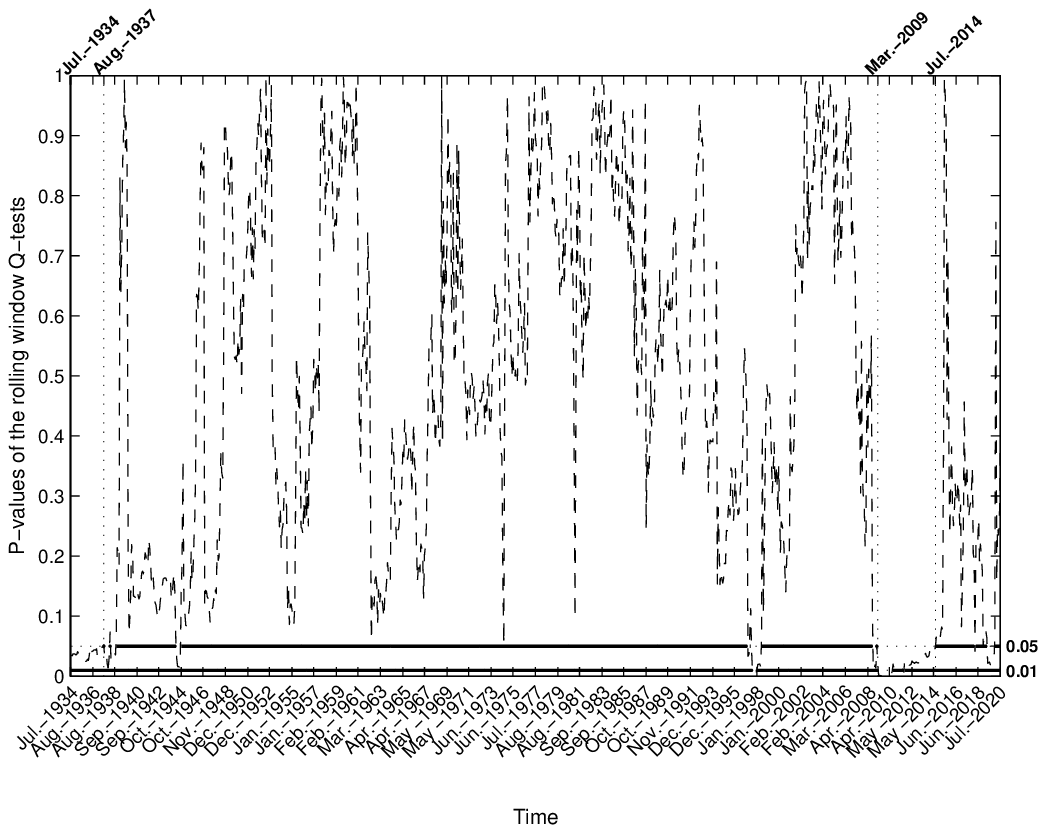}
\end{tabular}
\caption{The sample autocorrelation coefficient $\hat \rho_1$ computed by the moving window method for the monthly S\&P500 returns with $1$\% confidence bounds (left) and the $p$-values of the rolling window Q-test for residual autocorrelation (right).}
\label{fig:SP500}
\end{figure*}

\begin{figure*}[th!]
\begin{tabular}{cc}
\includegraphics[width = 0.5\textwidth]{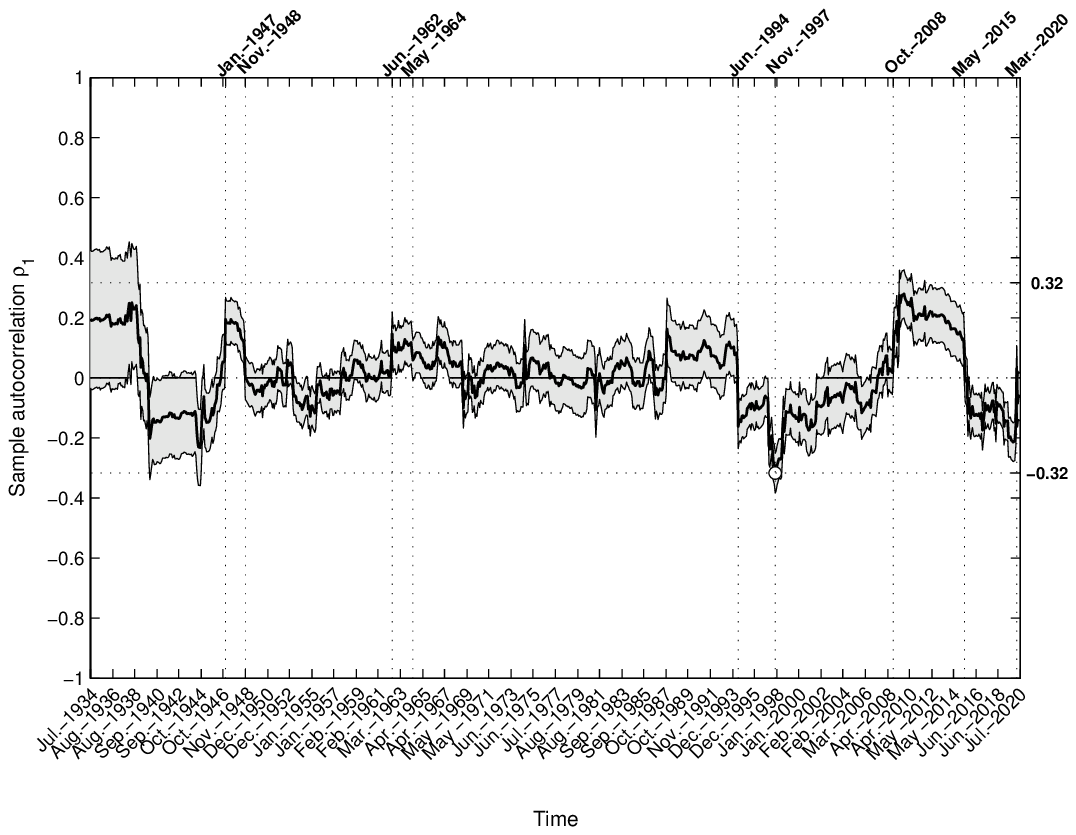} &
\includegraphics[width = 0.5\textwidth]{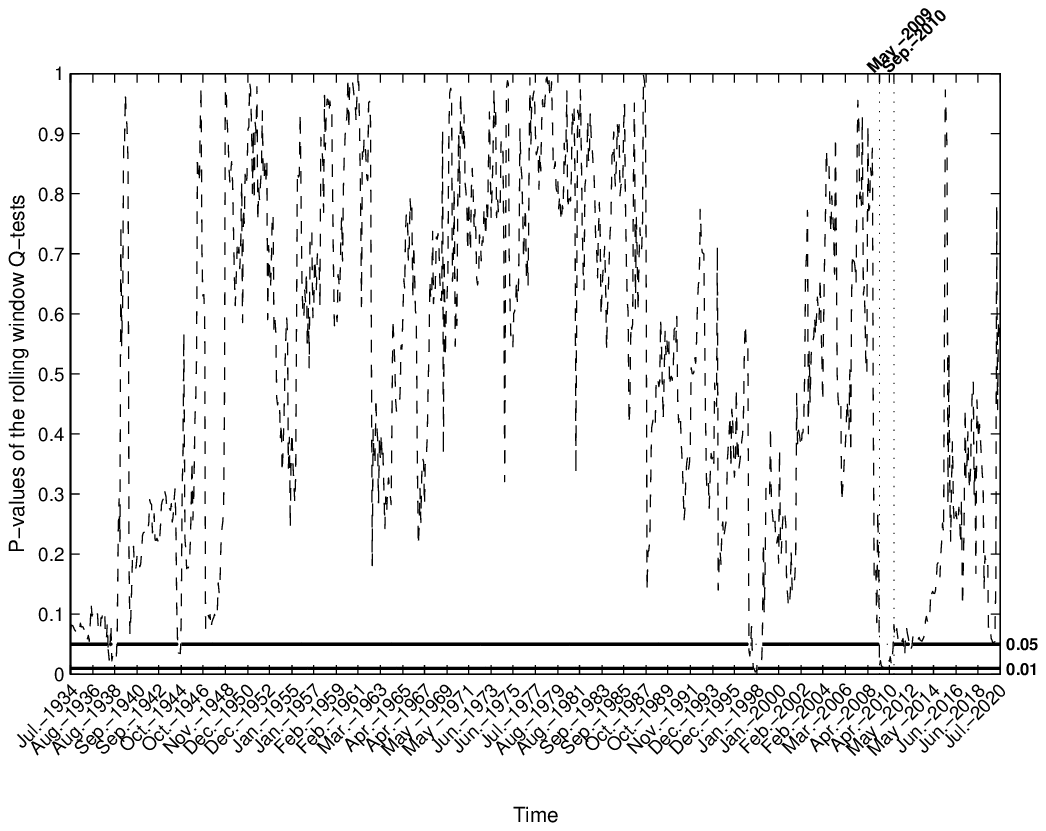}
\end{tabular}
\caption{The sample autocorrelation coefficient $\hat \rho_1$ computed by the moving window method for the monthly DJIA returns with $1$\% confidence bounds (left) and the $p$-values of the rolling window Q-test for residual autocorrelation (right).}
\label{fig:DJIA}
\end{figure*}

With the purpose of our research, we follow the moving window method for tracking the sample autocorrelation coefficient at lag $l=1$ as explained in~\cite{2004:Lo,2009:ItoSugiyama}. This provides insight into the underlying dynamics
of the time-varying AR(1) coefficient. Let $w$ denotes the window size where $w<N$ and $N$ is a number of observations. Define $N-w+1$ sub-samples of size $w$ by $\{y(t_{k-w+1}), \ldots, y(t_k) \}$ for $k = w, \ldots, N$ and then compute the first order autocorrelation to each sub sample. The window size is set to $w = 80$ and the corresponding $1$\% confidence bounds are calculated. The market is weak-form efficient when $\hat \rho_1$ equals zero (within the chosen confidence interval). To enhance any insight, we also provide the rolling window Ljung-Box Q-test results for residual autocorrelation at lag $l=1$ and, concurrently, plot the time-varying $p$-values. For $p$-values less than $\alpha$, say $\alpha = 0.01$ or $\alpha = 0.05$, the null hypothesis of uncorrelated residuals is rejected at the related significance level $\alpha$. This implies some degree of inefficiency exists during that time period.
Fig.~\ref{fig:SP500} illustrates a time varying structure of the $\hat \rho_1$ obtained for the S\&P500 stock index series together with the $p$-values of the rolling window Q-test for residual autocorrelation. Similarly, Fig.~\ref{fig:DJIA} summarizes the results obtained for the DJIA series.

The results of the empirical study yield a few important conclusions. First, we received a strong evidence of the time-varying nature of weak-form U.S. market efficiency. Indeed, the left-hand graphs of Figs.~\ref{fig:SP500} and~\ref{fig:DJIA} illustrate the sample autocorrelation coefficient $\hat \rho_1$ computed by the moving window method, which is a proxy for an evolving market efficiency level. The associated $1$\% confidence bounds are also presented. It is clearly seen that the estimated levels of market efficiency fluctuate around the critical level of zero except during the period related to the U.S. Great Depression and the 2008-2009 global financial crisis for both samples.  A short period of inefficiency is also observed during 1998s. We emphasize that a window size of $w=80$ months allows us to monitor the process from July 1934 till July 2020. During this period, the highest level of inefficiency for the S\&P500 index is observed on August 2009 with an estimated value of $\hat \rho_1 =  0.46$. Meanwhile, for the DJIA index we obtain an estimated value of $\hat \rho_1 =  -0.34$ observed on January 1998 and $\hat \rho_1 =  0.32$ on August 2009.

Secondly, it is clearly seen that the longest period of inefficiency since the beginning of our empirical study is detected from October 2008 till May 2015 for both the S\&P500 and DJIA indices. It corresponds well to an instability during the global 2008-2009 financial crisis. We also clearly observe the inefficiency period from July 1934 till August 1938 that is related to the U.S. Great Depression.

Thirdly, the efficiency levels recovered for both the S\&P500 and DJIA samples show similar patterns in their dynamics. Having analyzed the left-hand graphs of Figs.~\ref{fig:SP500} and~\ref{fig:DJIA}, we conclude that the U.S. market starts to show some signs of the recovering process toward weak-form efficiency regime  since August 1938. However, it was highly volatile during the Second World War and showed the period of inefficiency since May 1944 till October 1944. From the end of 1944 and till the end of 2008 the U.S. stock market is nearly always weak-form efficient at the $1$\% confidence level, although some short period of anomalies is detected in 1998s. From the estimated dynamics illustrated by the left-hand graphs of Figs.~\ref{fig:SP500} and~\ref{fig:DJIA}, we also conclude that the processes breached the efficiency levels in May 2009 for both samples and remained inefficient for nearly six years. This conclusion is substantiated for the S\&P500 sample by the right-hand graph of Fig.~\ref{fig:SP500} where the $p$-values of the related rolling window Ljung-Box Q-test for residual autocorrelation at lag $l=1$ are represented. This gives us an alternative version of the market efficiency investigation. It is interesting to note that the DJIA index showed the recovering process in a more rapid way than the S\&P500 index at the $5$\% confidence level. Indeed, following the right-hand graph of Fig.~\ref{fig:DJIA}, the estimated efficiency level of the DJIA index came back to an efficiency regime on September 2010 at the $5$\% confidence level, meanwhile the S\&P500 index remained inefficient until July 2014 at the same $5$\% confidence level.

To summarize the results of our empirical study, we conclude that the U.S. market is always weak-form efficient at the 5\% significance level, except two large  periods of inefficiency detected from July 1934 till August 1938 and from May 2009 till July 2014. Thus, we conclude that the U.S. Great Depression and the 2008-2009 global financial crisis had a dramatic impact on the U.S. stock market. It is also worth noting that our conclusion here concurs with the published results in~\cite{2009:ItoSugiyama,2011:Kim:AMH}, where different methods were used to investigate the AMH in the U.S. market. Finally, it is interesting to note that the most recent recession related to the first wave of the COVID-19 pandemic produced some significant movements in the level of market efficiency. More precisely, a movement toward inefficiency is observed on January 2020 for both samples under examination. Following Figs.~\ref{fig:SP500} and~\ref{fig:DJIA}, the U.S. market is highly volatile during the first half of 2020 with some swift returns to a weak-form efficiency regime. Our findings also reveal that the U.S. market is weak-form inefficient at the end of our empirical study on July 2020 at the 5\% significance level. This seems to be a consequence of the first wave of the COVID-19 pandemic.

The goal of this paper is to suggest the estimator appropriate for tracking the dynamics of the first order autoregression coefficient illustrated by the left-hand graphs of Figs.~\ref{fig:SP500} and~\ref{fig:DJIA} modeled by a time-varying AR(1) processes with both homoskedastic and heteroscedastic  conditional variance assumptions.

\section{Modeling and estimation of market efficiency dynamics} \label{Section:models}

To suggest an appropriate model along with estimation method, we come back to Table~\ref{tab:models} where the survey of various model specifications is provided. It becomes clear that a common technique  to modeling and tracking the {\it time-varying} level of market efficiency under the AMH
is to apply AR($n$) models with {\it time-dependent} coefficients to a return history. To design the related estimation method, we first start with the case of constant regression coefficients given by equation~\eqref{AR:n}. Following~\cite[Section~4.4.3]{2015:Grewal:book}, it can be represented in the form of linear predictive models as follows:
{\small
\begin{align*}
\begin{bmatrix}
y(t_{k+1})\\
y(t_{k})\\
y(t_{k-1})\\
\vdots \\
y(t_{k-n+2})
\end{bmatrix}
=
\begin{bmatrix}
\beta_1 & \beta_2 & \ldots & \beta_{n-1} & \beta_n \\
1 & 0 & \ldots & 0 & 0\\
0 & 1 & \ldots & 0 & 0 \\
\vdots & \vdots & \ddots & \vdots & \vdots \\
0 & 0 & \ldots & 1 & 0
\end{bmatrix}
\begin{bmatrix}
y(t_{k})\\
y(t_{k-1}) \\
y(t_{k-2})\\
\vdots \\
y(t_{k-n+1})
\end{bmatrix}
+
\begin{bmatrix}
\varepsilon(t_{k+1})\\
0 \\
0 \\
\vdots \\
0
\end{bmatrix}
\end{align*}
}
where $x(t_{k}) = [y(t_{k}), y(t_{k-1}), \ldots, y(t_{k-n+1})]^{\top}$ is the state vector that can be simulated by equation $x(t_{k+1}) = F x(t_{k})$ for any $t_k$ where $k = n, \ldots, N$ with the constant transition matrix $F$ defined above and the initial state $x(t_{0}) = [y(t_{n}), y(t_{n-1}), \ldots, y(t_{1})]^{\top}$. As can be seen, the transition matrix is parameterized by a constant vector of system parameters $\beta = [\beta_1, \ldots, \beta_n]$. Both the state vector and the system parameters might be estimated from a return series by representing the model in a linear state-space form and, next, by applying the maximum likelihood estimation (MLE) procedure along with the classical KF. Indeed, since the only first entry of the state vector is observed at each $t_k$, which is the return $y(t_k)$, the measurement equation is given as follows:
\begin{align}
z(t_{k}) =
\begin{bmatrix}
1 & 0 & \ldots & 0 & 0
\end{bmatrix}
x(t_{k}) + v(t_{k}), \; v(t_k) \sim {\cal N}(0,R) \label{reg:meas}
\end{align}
where the measurements are assumed to be noisy and one should define almost ``exact measurement" case because $z(t_k)$ is, in fact, our return $y(t_k)$ observed at time $t_k$. It is worth to mention that the setting $R=0$, which corresponds to the``exact measurement" scenario, yields to a failure of any KF-like estimation method due to its numerical instability. Hence, in our numerical experiments and empirical studies, we set $R=10^{-6}$ that allows us to simulate the almost ``exact measurement" case.

\subsection{Estimating time-varying autoregressive models} \label{Section:AR}

It is clearly seen that the AR($n$) process with time-varying coefficients $\beta(t_k) = [\beta_1(t_k), \ldots, \beta_n(t_k)]$, which are themselves modeled as separate stochastic processes, leads to a nonlinear discrete-time stochastic system. Despite this fact, the related models are often estimated by the classical KF and linear smoothing algorithms where formula~\eqref{AR:n:time} is considered as an observation equation of the state-space model representation. This means that the observation matrix is defined by $H_k = [y(t_{k}) | y(t_{k-1}) | \ldots | y(t_{k-n+1})]$. For instance, the readers are referred to the discussion in~\cite{1997:Emerson,1999:Zalewska,2009:ItoSugiyama} and some other papers. We stress that the inappropriate problem statement might diminish an estimation quality of the filter utilized for recovering the hidden state process $\beta(t_k)$ from a return history. Recall, if  the time-varying AR(1) model is utilized to test weak-form efficiency, then a degree of market efficiency is modeled by the first-order time-varying regression coefficient $\beta_1(t_k)$. Consequently, the hidden underlying processes illustrated by the left-hand graphs of Figs.~\ref{fig:SP500} and~\ref{fig:DJIA} should be recovered by the applied estimator from the returns available. In this section, we show how an accurate estimation method can be constructed within nonlinear Bayesian filtering framework.

Let us consider the AR(1) process in~\eqref{AR:n:time} with time-varying coefficient $\beta_1(t_k)$ following the random walk in~\eqref{Beta:RW}. Taking into account the discussion above, the model with homoskedastic conditional variance assumption can be written in the state-space form (the details can be found in Appendix) as follows:
\begin{align}
x_{k+1}&:=\begin{bmatrix}
\beta_1(t_{k+1}) \\ y(t_{k+1})
\end{bmatrix} = \left\{\begin{matrix}
\beta_1(t_{k}) + w(t_{k+1}) \\
\beta(t_{k+1})y(t_k) + \varepsilon(t_{k+1})
\end{matrix}\right. \label{eq:model1} \\
& =
\left\{\begin{matrix}
\beta_1(t_{k}) + w(t_{k+1}) \\
\beta_1(t_{k})y(t_k) + y(t_k) w(t_{k+1}) + \varepsilon(t_{k+1})
\end{matrix} \right. \nonumber \\
& = \underbrace{\begin{bmatrix}
\beta_1(t_{k}) \\
\beta_1(t_{k})y(t_k)
\end{bmatrix}}_{f(x_k,t_{k+1},t_k)} +  \underbrace{\begin{bmatrix}
1 & 0 \\
y(t_k)&  1
\end{bmatrix}}_{G_{k}}
\underbrace{\begin{bmatrix}
w(t_{k+1}) \\
\varepsilon(t_{k+1})
\end{bmatrix}}_{u_{k+1}}, \quad
 \underbrace{\begin{bmatrix}
\sigma_{w}^2 & 0 \\
0 &  \sigma_{\varepsilon}^2
\end{bmatrix}}_{Q}  \nonumber
\end{align}
where the entries of the state vector $x_k = [x_1(t_k), x_2(t_k)]^{\top}$ are defined by $x_1(t_k):=\beta_1(t_k)$ and $x_2(t_k):=y(t_k)$, i.e. $x_k = [\beta_1(t_k), y(t_k)]^{\top}$. The value $y(t_k)$ is the mean-adjusted log return at $t_k$, $k=1,\ldots, N$. The initial state is then defined by $\bar x_0 = [\beta_1^{(0)}, y(t_1)]^{\top}$.  Meanwhile, the first entry $\beta_1^{(0)}$ is a proxy for the autocorrelation coefficient $\rho_1$ and, hence, its absolute value can not exceed one. It also corresponds to the locally stationary process in~\eqref{AR:n:time} as discussed in~\cite{1997:Dahlhaus}. Following~\cite{1984:Harvey}, we consider the initial efficiency level $\beta_1^{(0)}$ as an extra unknown system parameter that should be estimated from the data while calibrating the model by the MLE procedure. We also set the initial error covariance matrix by $\Pi_0 = I_2$ and refer to the discussion presented in Appendix, which is related to the choice of the filter initials and parameter estimation technique.

\begin{figure}
\includegraphics[width=0.5\textwidth]{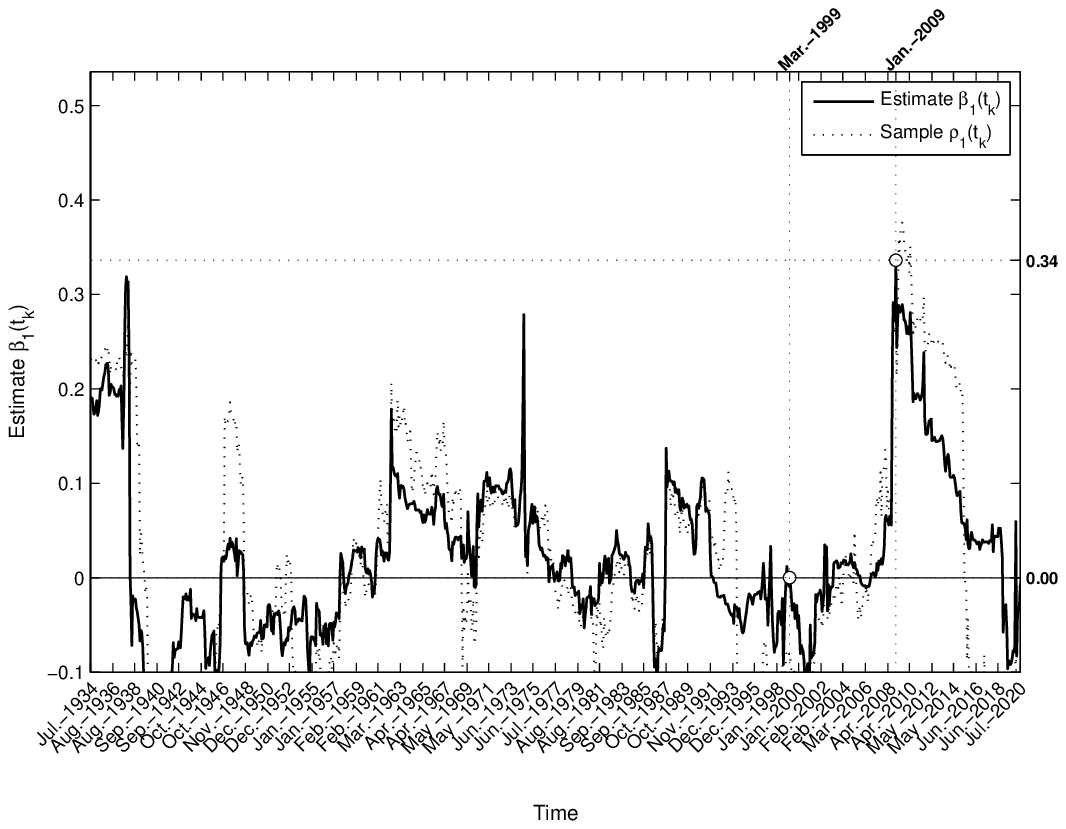}
\caption{The estimated trajectory $\hat \beta_1(t_k)$ modeled by~\eqref{eq:model1} and tracked by the EKF and the sample autocorrelation coefficient $\hat \rho_1$ computed by the moving window method for the monthly S\&P500 returns.} \label{beta1:AR1RW:graph}
\end{figure}

The measurement equation in the state-space model (see also equation~\eqref{eq:state:space2} in Appendix) is linear and has the form:
\begin{align}
z(t_{k}) = \begin{bmatrix}
0 &  1
\end{bmatrix}
\begin{bmatrix}
x_1(t_k) \\
x_2(t_k)
\end{bmatrix}
+ v(t_{k}), \; v(t_k) \sim {\cal N}(0,R), \label{AR1RW:measurement}
\end{align}
 where only the second entry of the state vector is observed; see formula~\eqref{reg:meas} and the discussion below that formula concerning the noise covariance matrix $R$.

When the model is casted in the state-space form, the unknown state vector can be estimated by nonlinear Bayesian filters. The extended Kalman filter (EKF) is a simple and attractive nonlinear estimator, which provides a good balance between accuracy and computational cost, especially for the state-space models with linear measurement equation. The EKF equations are summarized by formulas~\eqref{state_new}~-- \eqref{kf:f:X} where the nonlinear drift functions and Jacobian matrix are defined by
\begin{align}
f\bigl(x_k,t_{k},t_{k+1}\bigr) & =
\begin{bmatrix}
x_{1}(t_{k}) \\
x_1(t_{k}) x_2(t_{k})
\end{bmatrix}, & \!\!\!
F_k & =
\begin{bmatrix}
1 & 0\\
x_2(t_{k}) & x_1(t_{k})
\end{bmatrix}. \label{EKF:model:1}
\end{align}

Finally, the system parameters $\theta = [\sigma_{w_1}^2, \sigma_{\varepsilon}^2,\beta_1^{(0)}]$ are unknown and should be estimated from the real data available. In other words, the model should be fitted to the data and then the EKF is applied for recovering the hidden state vector. In practice, these two problems can be solved simultaneously by applying the {\it adaptive} EKF as explained in Appendix. Thus, the unknown system parameters are estimated by the MLE and the hidden state vector is recovered by the EKF.

 \begin{figure}
\includegraphics[width=0.5\textwidth]{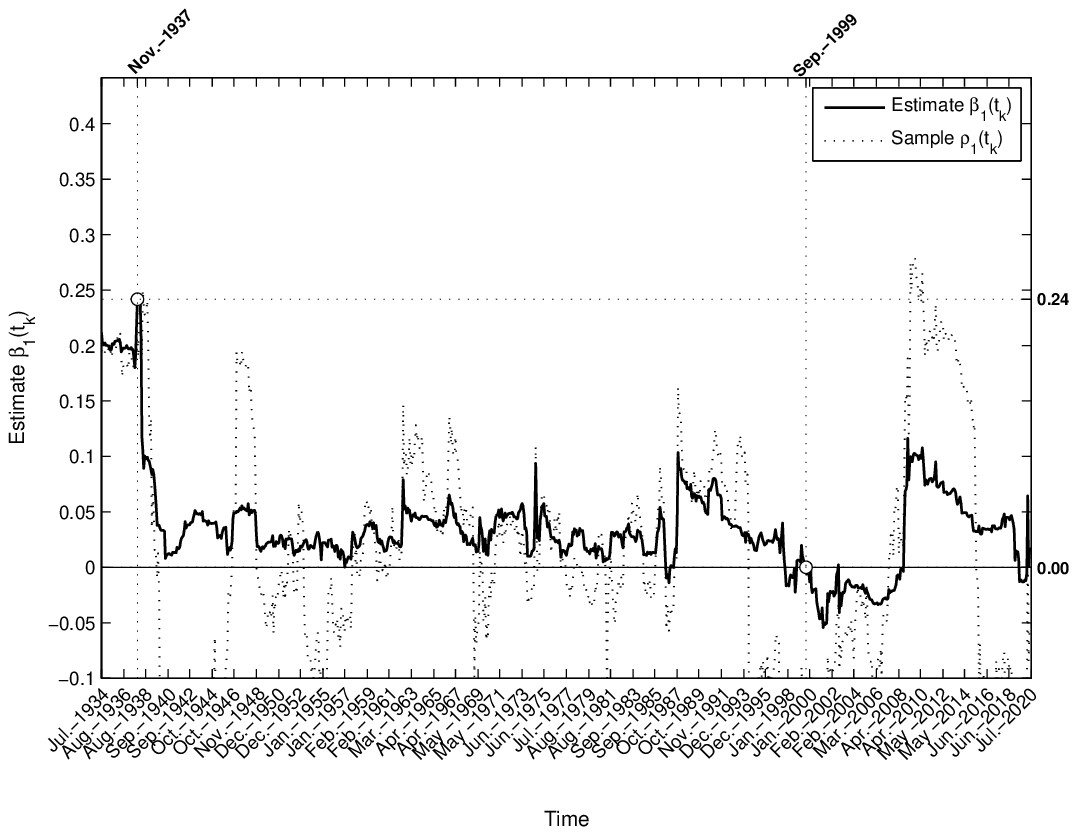}
\caption{The estimated trajectory $\hat \beta_1(t_k)$ modeled by~\eqref{eq:model1} and tracked by the EKF and the sample autocorrelation coefficient $\hat \rho_1$ computed by the moving window method for the monthly DJIA returns.} \label{beta1:AR1RW:graph:DJIA}
\end{figure}

\begin{table*}[ht!]
\caption{Estimation results for monthly S\&P500 and DJIA mean-adjusted returns with standard errors given in parentheses. The lowest AIC indicates the preferred model by $\blacktriangledown$ estimated by the novel EKF-based method, meanwhile $\blacktriangle$ denotes the best fitted model within the classical KF approach.} \label{Table:estimate:AR1RW}                                                                                                                                                                                                                                                                                                                                                                                                                                                                                                                                                                                                                                                                {\small
\begin{tabular}{l||c|c||c|c||c|c||c|c}
\toprule
{\bf Model} & \multicolumn{4}{c||}{\bf S\&P500 index: AR(1) with time-varying coefficients} & \multicolumn{4}{c}{\bf DJIA index: AR(1) with time-varying coefficients}  \\
\cline{2-9}
& \multicolumn{2}{c|}{({\it without trend})} & \multicolumn{2}{c||}{({\it with trend})} & \multicolumn{2}{c|}{({\it without trend})} & \multicolumn{2}{c}{({\it with trend})}  \\
& \multicolumn{2}{c|}{given by eqs.~\eqref{eq:model1}} & \multicolumn{2}{c||}{given by eqs.~\eqref{eq:model2}} & \multicolumn{2}{c|}{given by eqs.~\eqref{eq:model1}} & \multicolumn{2}{c}{given by eqs.~\eqref{eq:model2}}  \\
\cline{2-9}
{\bf Method} & novel EKF & KF in~\cite{2009:ItoSugiyama} &  novel EKF & KF in~\cite{2009:ItoSugiyama}  & novel EKF & KF in~\cite{2009:ItoSugiyama} & novel EKF & KF in~\cite{2009:ItoSugiyama} \\
\toprule
$\mu_{\beta_1}$           &  ---    & ---      & -0.00022       & -0.00022 &   ---   &   --- & -0.00014 & -0.00014	 \\
                          &         &          & ({\it 0.00062})  &  ({\it 0.00063}) & & & ({\it  0.00008}) &  ({\it  0.00008}) \\
$\sigma_{\varepsilon}^2$  & 0.00284 &  0.00284 & 0.00284     &   0.00284  &  0.00275  & 0.00275  & 0.00277 & 0.00276 \\
                          & ({\it  0.00012}) & ({\it  0.00012}) & ({\it 0.00012}) & ({\it 0.00012}) & ({\it  0.00012}) & ({\it  0.00012}) & ({\it  0.00012}) & ({\it  0.00012}) \\
$\sigma_{w}^2$            & 0.00040         & 0.00041   &  0.00039 & 0.00040   & 0.00007 & 0.00007 & 0.00001 & 0.00001  \\
                          & ({\it 0.00039}) & ({\it 0.00039}) & ({\it 0.00037}) & ({\it 0.00040}) & ({\it 0.00016}) & ({\it 0.00016}) & ({\it 0.00001}) & ({\it 0.00001}) \\
$\beta_1^{(0)}$           & 0.22299         &  0.22236 &   0.23032 &	0.23012	 & 0.15342 & 0.15349 & 0.12582 & 0.12479 \\
 ({\it initial})          & ({\it 1.07033}) &  ({\it 0.32182})  & ({\it 0.94585}) &  ({\it 0.94585}) & ({\it 0.00001}) & ({\it 0.00001}) & ({\it 0.41257}) & ({\it 0.00001}) \\
\hline
$\max {\cal L}$ & 2687.96  & 2690.88 & 2688.02 & 2690.94 & 2710.50 & 2713.41 & 2711.06 & 2713.95 \\
AIC             & -5369.93$\blacktriangledown$ & -5375.76$\blacktriangle$  &  -5368.06 &  -5373.88 & -5415.01$\blacktriangledown$  & -5420.82$\blacktriangle$ & -5414.12 & -5419.91 \\
\hline
\hline
$\| \hat \rho_1(t) - \hat \beta_1(t)\|_{\infty}$ & 0.2967 &   0.3026	 &  0.2975 & 0.3033	& 0.2985 & 0.2999 & 0.3363	& 0.3364 \\
\bottomrule
\end{tabular}
}
\end{table*}

The estimates calculated for both the S\&P500 and DJIA series are summarized in Table~\ref{Table:estimate:AR1RW}. The estimated paths $\hat \beta_1(t_k)$ recovered by the novel EKF-based approach are illustrated by Figs.~\ref{beta1:AR1RW:graph} and~\ref{beta1:AR1RW:graph:DJIA}, respectively. It is clearly seen that the trajectory $\left\{ \hat \beta_1(t_k)\right\}$ recovered by the filtering method adequately follows the sample autocorrelation coefficient $\left\{\hat \rho_1 (t_k)\right\}$ computed by the moving window method in Section~\ref{sec:data}, which is a proxy for an evolving market efficiency level. To substantiate an estimation quality of the novel EKF-based estimation approach, we additionally provide a comparative study of the new method and the classical KF framework previously applied in~\cite{2009:ItoSugiyama} and other papers. To assess a quality of two estimation methods, we calculate the norm $\| \hat \rho_1(t) - \hat \beta_1(t)\|_{\infty}$, which gives some insights into the market efficiency level estimation process. In summary, these values allow us to decide about the accuracy of each estimation method applied and answer the question on how adequate did they recover the sample autocorrelation coefficient. Finally, we examine two time-varying AR(1) models for evolving weak-form market efficiency tests.  The estimation errors are computed for both filtering methods applied for each model under examination and the obtained results are summarized in Table~\ref{Table:estimate:AR1RW}.

Let us explore a comparative study of the estimation methods, first. For a fair comparison, each adaptive filtering
strategy utilizes precisely the same sampled data, i.e. the S\&P500 and the DJIA  monthly returns, the same
initial system parameters $\theta^{(0)} = [\sigma_{w_1}^2, \sigma_{\varepsilon}^2,\beta_1^{(0)}] := [0.01,0.1,0]$ for the optimization method implemented and the same optimization code, which is taken to be the MATLAB built-in function {\tt fmincon}. Having analyzed the results collected in the left and
right panels of Table~\ref{Table:estimate:AR1RW}, we see that all the methods under examination
recover the hidden state vector of each model equally well, i.e. with similar accuracies. Besides, the maximum reconstruction errors are small for all models estimated and for both the S\&P500 and the DJIA samples. More importantly, the reconstruction errors of the novel EKF-based estimation method are smaller than for the classical KF procedure previously used in econometric literature. This conclusion holds true for both models examined and both return series explored. For instance, the reconstruction error of the novel EKF approach while estimating model~\eqref{eq:model1} for the S\&P500 index has an estimated value of $\approx 0.29$, which is smaller than $\approx 0.30$ obtained by the KF framework for the same model and the same return series. Having analyzed the values calculated for the DJIA series, we obtain the same conclusion about a higher estimation quality of the novel EKF-based method over the standard KF technique. Finally, we need to emphasize that although the EKF approach outperforms the classical KF framework, the difference in their estimation quality is small. It should be stressed that the simple time-varying AR(1) models with homogeneous conditional variance assumption are investigated in this Section. We expect to obtain a more significance difference in estimation quality of these two filtering approaches while exploring more sophisticated tests based on time-varying AR(1) models with heteroscedastic conditional variance assumption.

Let us next examine a choice of the models to perform the tests for evolving market efficiency. It is interesting to note that the U.S. market has been explored under the AMH in~\cite{2009:ItoSugiyama}. The same time-varying AR(1) models with homogeneous conditional variance assumption were utilized and estimated by the classical Kalman filter. It was found out that the relative degree of the U.S. market efficiency varies through time without trend from January 1955 till February 2006 as reported in the cited paper. To decide whether or not the same conclusion holds for the S\&P500 and the DJIA returns for the period from November 1927 to July 2020 examined in this paper, we additionally explore the time-varying AR(1) with a constant trend. Instead of the first dynamic equation in~\eqref{eq:model1}, we have
\begin{align}
\beta_1(t_{k+1}) & = \beta_{1}(t_{k}) + \mu_{\beta_1} + w_{1}(t_{k+1}), \; w_{1}(t_{k+1})\sim {\cal N}(0,\sigma_{w}^2) \label{AR1RW:state1:trend}
\end{align}
where $\mu_{\beta_1}$ is an extra system parameter that should be estimated.

In summary, we have the following state-space model:
\begin{align}
x_{k+1}&:=\begin{bmatrix}
\beta_1(t_{k+1}) \\ y(t_{k+1})
\end{bmatrix} = \left\{\begin{matrix}
\beta_1(t_{k}) + \mu_{\beta_1} + w(t_{k+1}) \\
\beta(t_{k+1})y(t_k) + \varepsilon(t_{k+1})
\end{matrix}\right. \label{eq:model2} \\
& =
\left\{\begin{matrix}
\beta_1(t_{k}) + \mu_{\beta_1} + w(t_{k+1}) \\
\left(\beta_1(t_{k})+ \mu_{\beta_1}\right) y(t_k) + y(t_k) w(t_{k+1}) + \varepsilon(t_{k+1})
\end{matrix} \right. \nonumber \\
& = \underbrace{\begin{bmatrix}
\beta_1(t_{k}) + \mu_{\beta_1} \\
\beta_1(t_{k})y(t_k) + \mu_{\beta_1}y(t_k)
\end{bmatrix}}_{f(x_k,t_{k+1},t_k)} +  \underbrace{\begin{bmatrix}
1 & 0 \\
y(t_k)&  1
\end{bmatrix}}_{G_{k+1}}
\underbrace{\begin{bmatrix}
w(t_{k+1}) \\
\varepsilon(t_{k+1})
\end{bmatrix}}_{u_{k+1}},
 \underbrace{\begin{bmatrix}
\sigma_{w}^2 & 0 \\
0 &  \sigma_{\varepsilon}^2
\end{bmatrix}}_{Q}  \nonumber
\end{align}
where $\theta = [\sigma_{w_1}^2, \sigma_{\varepsilon}^2,\beta_1^{(0)},\mu_{\beta_1}]$ is the unknown system parameter vector that should be estimated together with the unknown state.

The model summarized by formulas~\eqref{eq:model2} is estimated in the same way as discussed above where instead of~\eqref{EKF:model:1} we have
\begin{align}
& f\bigl(x_k,t_{k},t_{k+1}\bigr)  =
\begin{bmatrix}
x_{1}(t_{k}) + \mu_{\beta_1}\\
x_1(t_{k}) x_2(t_{k}) + \mu_{\beta_1}x_2(t_{k})
\end{bmatrix}, \label{EKF:model2:1} \\
& F_k  =
\begin{bmatrix}
1 & 0\\
x_2(t_{k}) & x_1(t_{k}) + \mu_{\beta_1}
\end{bmatrix}, \quad G_k =
\begin{bmatrix}
1 & 0\\
x_2(t_{k}) & 1
\end{bmatrix}.
\end{align}

The obtained estimates are summarized in Table~\ref{Table:estimate:AR1RW} as well. We use the Akaike information criterion (AIC) to choose the best fitted model that turns out to be the model summarized in~\eqref{eq:model1}. To ensure a robust comparison, all models are estimated using the same filtering method (the novel EKF and the classic KF), the same parameter estimation strategy (the method of maximum likelihood), the same optimization method (the built-in optimizer {\tt fmincon} in MATLAB), the same initial parameter values, and the same initialization of the filter. Both estimation methods indicate model~\eqref{eq:model1} as the best fitted choice. This means that the S\&P500 and the DJIA efficiency levels travel through time without trend over the period under examination, i.e. from November 1927 to July 2020. This conclusion is sustainable, i.e. it does not depend on the filtering method applied to calibrate the models.

Additionally,  we provide the likelihood ratio tests for the hypothesis $H_0: \mu_{\beta_1} = 0$ tested against $H_a: \mu_{\beta_1} \ne 0$ at the $1\%$ confidence level to justify the conclusion concerning a possible trend in the time-varying market efficiency level. The result of the likelihood ratio test indicates that the hypothesis $H_0: \mu_{\beta_1} = 0$  is not rejected at a $1\%$ confidence level, i.e. the restricted model with $\mu_{\beta_1} = 0$ fits the data better compared to its counterparts with $\mu_{\beta_1} \ne 0$, i.e. model~\eqref{eq:model1} is preferable for practical use since it contains fewer parameters. In other words, we substantiate the previous finding for the S\&P500 index that the relative degree of U.S. market efficiency varies through time without trend over the period under examination, i.e. from November 1927 to July 2020. We also discover that the efficiency level of the DJIA index travels through time without trend as well. Finally, Figs.~\ref{beta1:AR1RW:graph} and~\ref{beta1:AR1RW:graph:DJIA} illustrate the $\hat \beta_1(t_k)$, $k=1,\ldots,N$, recovered by the best fitted model from the S\&P500 and the DJIA monthly return series, respectively, against the sample autocorrelation coefficients $\hat \rho(t_k)$ calculated for the same return series by the moving window procedure.

Having analyzed Figs.~\ref{beta1:AR1RW:graph} and~\ref{beta1:AR1RW:graph:DJIA}, we clearly see that the estimated trajectories of $\hat \beta_1(t_k)$, $k=1,\ldots,N$, are quite similar to those of autocorrelations $\hat \rho_1$, $k=w,\ldots,N$. This additionally substantiates a good estimation quality of the novel EKF-based estimation framework proposed in this paper. In summary, the suggested model and the related filtering method are capable for adequate tracking the value that we are looking for. Nevertheless, there is a space for further improvements in terms of developing more sophisticated models and estimation methods for monitoring the time-varying market efficiency and detecting the periods of anomalies.

\subsection{Nonlinear serial dependence and GARCH process}

In previous section, time-varying AR(1) models with the homoskedastic conditional variance assumption have been investigated to perform the evolving market efficiency tests. Assumption about constant variance $\sigma_{\varepsilon}^2$ significantly restricts a capacity of modeling approaches for estimating the market efficiency level. It is widely acknowledge in the econometric and financial literature that the nonlinear serial dependence might be observed due to a dynamics of conditional variance process. More precisely, it is well known that an uncorrelated time series can still be serially dependent when they exhibit conditional heteroscedasticity or autocorrelation in the squared series. The models that allow to capture a time-varying nature of the volatility process (as well as the `clustering effect') are known as Autoregressive Conditional Heteroscedasticity (ARCH) originally defined in~\cite{1982:Engle} and later extended to generalized ARCH (GARCH) in~\cite{1986:Bollerslev}. As a result, the tests for weak-form  market efficiency with the implied conditional heteroscedasticity consist of AR($n$) model~\eqref{AR:n:time} with time-varying regression coefficients $\beta_{i,t}$, $i=1,\ldots,n$ modeled by random walk in~\eqref{Beta:RW}, and time-varying volatility process $\sigma_t^2:=\sigma_{\varepsilon}^2(t)$ followed the GARCH(p,q) specification:
\begin{align}
\sigma^2_t & =  \omega + \sum \limits_{l=1}^{p} a_l \varepsilon_{t-l}^2 + \sum \limits_{j=1}^{q} b_j \sigma^2_{t-j}. \label{garch:pq}
\end{align}

The test for weak-form market efficiency that takes into account heteroscedasticity of return series still yields the following condition: if $\beta_{1,t}=\beta_{2,t}= \ldots =\beta_{n,t}=0$, then equation~\eqref{AR:n:time} collapses to $y_t = \varepsilon_t$ where $\varepsilon_t \sim {\cal N}(0,\sigma^2_t)$ is a white noise. This implies that the log returns are serially uncorrelated in time but autocorrelation might be exhibited in the squared returns.

To build an appropriate model for evolving market efficiency, let us consider a general return series equation given by
\begin{align}
y_t & = \mu_t + \sigma_t\varepsilon_t, & \varepsilon_t & \sim {\cal N}(0,1) \label{return:eq}
\end{align}
where $\varepsilon_t \sim {\cal N}(0,1)$ is a white Gaussian noise process, $\mu_t$ is the conditional mean. The conditional variance $\sigma_t^2$ follows ARCH/GARCH process. With a purpose of testing time-varying market efficiency, the mean $\mu_t$ is modeled by time-varying AR($n$) given by formulas~\eqref{AR:n:time}, \eqref{Beta:RW} and, next, we restrict ourselves to  GARCH(1,1) process for the conditional variance equation.  Thus, we get the model of the following form:
\begin{align}
y_t & =  \beta_{1,t} y_{t-1} + \sigma_t \varepsilon_t, \quad \varepsilon_t \sim {\cal N}(0,1) \label{AR:n:time1} \\
\beta_{1,t} & = \beta_{1,t-1} + w_{1,t}, \quad  w_{1,t} \sim {\cal N}(0,\sigma_{w_1}^2) \label{Beta:RW1} \\
\sigma^2_t & = \omega + a_1 (y_{t-1} - \beta_{1,t-1} y_{t-2})^2 + b_1 \sigma^2_{t-1} \nonumber \\
& = \omega + a_1 (\sigma_{t-1}\varepsilon_{t-1})^2 + b_1 \sigma^2_{t-1} \label{garch:11}
\end{align}
where $\omega>0$, $a_1 \ge 0$, $b_1 \ge 0$ and $a_1 + b_1 <1$ ensures that the GARCH(1,1) process is stationary. The unconditional expectation of the conditional variance is constant and finite, and given by the formula $\omega/(1-a_1-b_1)$. The test for a weak-form market efficiency  implies the condition $\beta_{1,t}=0$ that yields the reduced formula $y_t = \sigma_t\varepsilon_t$ where $\varepsilon_t \sim {\cal N}(0,1)$ is a white Gaussian noise, i.e. the returns are serially uncorrelated in time.

Having denoted $x_1(t_k):=\beta_1(t_k)$, $x_2(t_k):=y(t_k)$ and $x_3(t_k):=\sigma^2(t_k)$, we define the state vector as follows: $x_k = [\beta_1(t_k),y(t_k),\sigma^2(t_k)]^{\top}$, and then express equations~\eqref{AR:n:time1}~-- \eqref{garch:11} in terms of $x_k =[x_1(t_k),x_2(t_k),x_3(t_k)]^{\top}$. When the model is written down in the state-space form~\eqref{eq:state:space2}, \eqref{eq:state:space3}, the adaptive EKF  procedure in~\eqref{state_new}~-- \eqref{kf:f:X} can be applied for the hidden state and system parameters estimation.

The measurement equation again has a linear form, i.e.
\begin{align}
z(t_{k}) = \begin{bmatrix}
0 &  1 & 0
\end{bmatrix}
\begin{bmatrix}
x_1(t_k) \\
x_2(t_k) \\
x_3(t_k)
\end{bmatrix}
+ v(t_{k}), \; v(t_k) \sim {\cal N}(0,R), \label{AR1RW:measurement}
\end{align}
because only the second entry of the state vector is observed; see formula~\eqref{reg:meas} and the discussion below that formula concerning the noise covariance matrix $R$.

\begin{figure}
\includegraphics[width=0.5\textwidth]{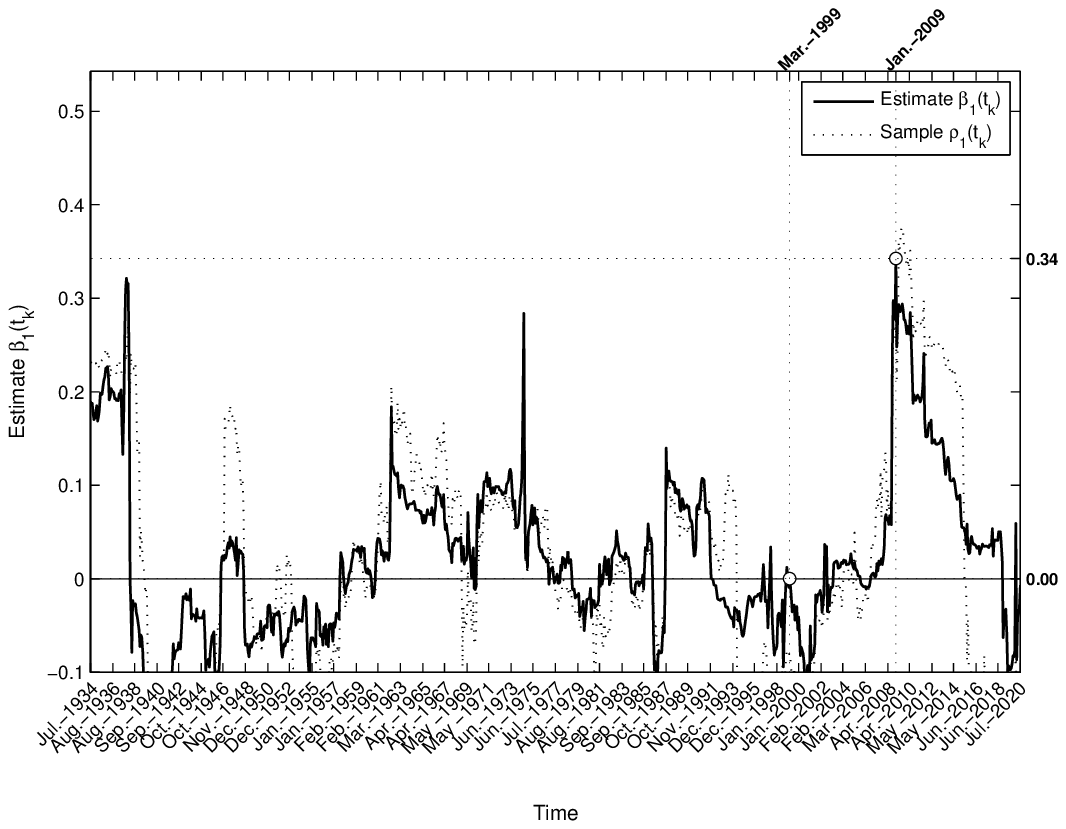}
\caption{The estimated trajectory $\hat \beta_1(t_k)$ modeled by eqs.~\eqref{AR:n:time1}~-- \eqref{garch:11}  and tracked by the EKF and the sample autocorrelation coefficient $\hat \rho_1$ computed by the moving window method for the monthly S\&P500 returns.} \label{beta1:General:graph}
\end{figure}

Unfortunately, the dynamic equation now has a more sophisticated form compared to the time-varying AR models with the assumed constant variance examined in Section~\ref{Section:AR}. Indeed, we have a complicated functional dependence at the process disturbance that yields a general equation~\eqref{eq:state:space3} instead of simple case in~\eqref{eq:state:space1}. We also note that the right-hand function in the dynamic equation~\eqref{eq:state:space3} depends on the current $x_k:=x(t_k)$, but may not depend on the future state values $x_{k+1}:=x(t_{k+1})$. Hence, equation~\eqref{AR1RW:state2:new} is transformed as required, i.e. we get
\begin{align}
x_1(t_{k+1})  & = x_{1}(t_{k}) + w_{1}(t_{k+1}), \; w_{1}(t_{k+1}) \sim {\cal N}(0,\sigma_{w_1}^2), \label{AR1RW:state1:new}\\
x_2(t_{k+1}) & = x_1(t_{k+1}) x_2(t_{k}) + \!\!\sqrt{x_3(t_{k+1})}\varepsilon(t_{k+1}) \nonumber \\
&= x_{1}(t_{k})x_2(t_{k}) + x_2(t_{k})w_{1}(t_{k+1}) \nonumber \\
& + \sqrt{\omega + a_1 x_3(t_{k}) \varepsilon^2(t_{k}) + b_1 x_3(t_{k})}\varepsilon(t_{k+1}), \label{AR1RW:state2:new} \\
x_3(t_{k+1}) & = \omega + b_1 x_3(t_{k}) + a_1 x_3(t_{k}) \varepsilon^2(t_{k}) \label{garch:11:new}
\end{align}
where $\varepsilon(t_{k+1}) \sim  {\cal N}(0,1)$ is a white noise and, hence, it is uncorrelated in time. It is also uncorrelated with $w_{1}(t_{k+1})$ and with the initial state $\bar x_0$. Below, we set $\bar x_0 = [\beta_1^{(0)}, y(t_1), \omega/(1-a_1-b_1)]^{\top}$ and $P_0 = I_3$ where $\theta = [\beta_1^{(0)},\omega,a_1,b_1,\sigma_{w}^2]$ is the unknown vector of system parameters that should be estimated from the return series.

To extract the hidden state process from measurements available, the EKF strategy presumes a linear approximation of the right-hand side function
$f(x_k,t_{k+1},t_k,u_{k+1})$ at the filtering state expectation
$\hat x_{k|k}$ of the hidden state $x_{k}$ and around $\E{u_{k+1}}=0$. Thus, we write down the linearized model~\eqref{eq2.2} as follows:
\begin{align*}
\tilde x_{k} & =\begin{bmatrix}
\hat x_{1,k|k} \\
\hat x_{1,k|k} \cdot \hat x_{2,k|k} \\
\omega + b_1 \hat x_{3,k|k}
\end{bmatrix} \!\!
+ \!\!\underbrace{\begin{bmatrix}
1 & 0 & 0\\
\hat x_{2,k|k} & \hat x_{1,k|k} & 0 \\
0 & 0 & b_1
\end{bmatrix}}_{\left.\partial f\bigl(x_{k},t_{k+1},t_k,0\bigr)/\partial x_k\right|_{\hat x_{k|k}}} \!\!\! \left(\begin{bmatrix}
x_1(t_k) \\
x_2(t_k) \\
x_3(t_k)
\end{bmatrix} \!\!
- \!\! \begin{bmatrix}
\hat x_{1,k|k} \\
\hat x_{2,k|k} \\
\hat x_{3,k|k}
\end{bmatrix}\right) \!\!\!\\
& + \underbrace{\begin{bmatrix}
1 & 0 & 0\\
\hat x_{2,k|k} & \sqrt{\omega + b_1 \hat x_{3,k|k}} & 0 \\
0 & 0 & 0
\end{bmatrix}}_{\left.\partial f\bigl(\hat x_{k|k},t_{k+1},t_k,u_{k+1}\bigr)/\partial u_{k+1}\right|_{0}}
\begin{bmatrix}
w_{1}(t_{k+1}) \\
\varepsilon(t_{k+1}) \\
\varepsilon(t_{k})
\end{bmatrix}
\end{align*}
where the noise covariance is $Q={\rm diag}\{\sigma_{w}^2,1,1\}$. To summarize, the Jacobian matrices used above are the following ones:
\begin{align*}
F_k =
\begin{bmatrix}
1 & 0 & 0\\
x_2(t_{k}) & x_1(t_{k}) & 0 \\
0 & 0 & b_1
\end{bmatrix},
G_k  = \begin{bmatrix}
1 & 0 & 0\\
x_2(t_k) & \sqrt{\omega + b_1 x_3(t_{k})} & 0 \\
0 & 0 & 0
\end{bmatrix}
\end{align*}
and should be computed at $\hat x_{k|k}$ in each filtering step.

\begin{figure}
\includegraphics[width=0.5\textwidth]{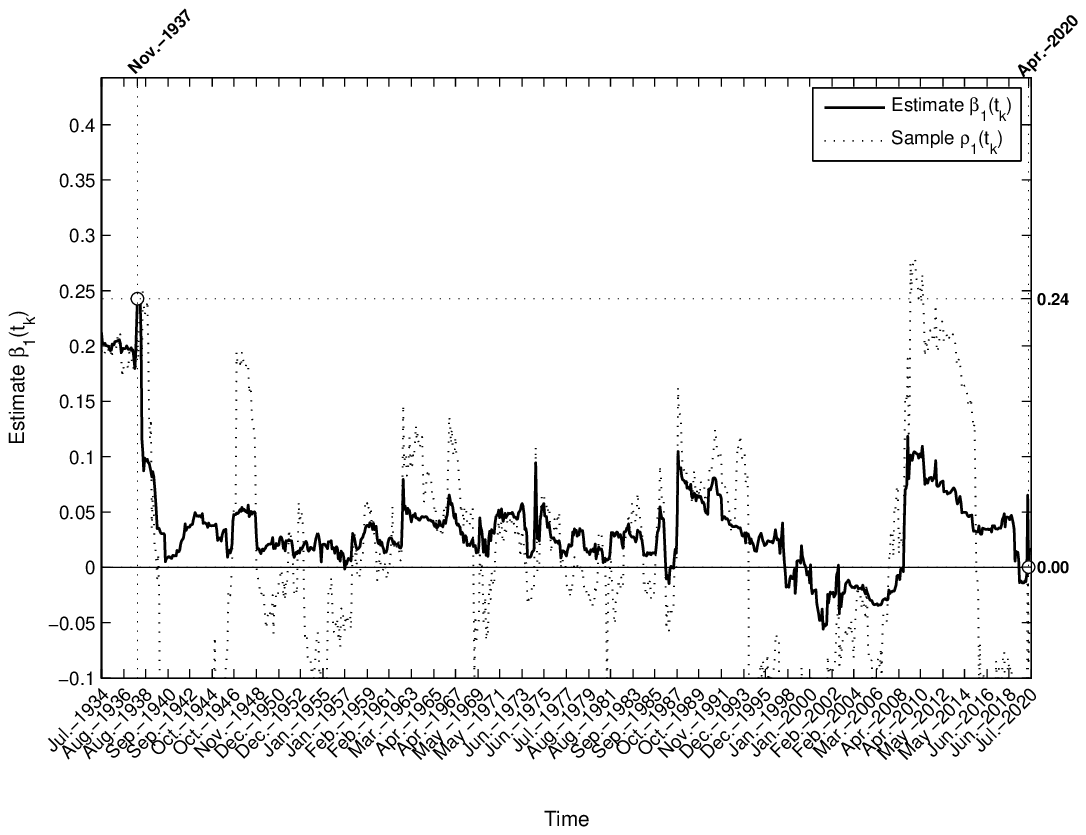}
\caption{The estimated trajectory $\hat \beta_1(t_k)$ modeled by eqs.~\eqref{AR:n:time1}~-- \eqref{garch:11}  and tracked by the EKF and the sample autocorrelation coefficient $\hat \rho_1$ computed by the moving window method for the monthly DJIA returns.} \label{beta1:General:graph:DJIA}
\end{figure}

\begin{table*}[ht!]
\caption{Estimation results of the evolving market efficiency test modeled by eqs.~\eqref{AR:n:time1}~-- \eqref{garch:11} with standard errors given in parentheses.} \label{Table:estimate:AR1RW:GARCH}                                                                                                                                                                                                                                                                                                                                                                                                                                                                                                                                                                                                                                                                \begin{tabular}{l||c|c||c|c}
\toprule
{\bf Series} & \multicolumn{2}{c||}{\bf S\&P500 index} & \multicolumn{2}{c}{\bf DJIA index} \\
\cline{2-5}
{\bf Method} & novel EKF & classical KF in~\cite{2008:Posta,1999:Zalewska} &  novel EKF & classical KF in~\cite{2008:Posta,1999:Zalewska} \\
\toprule
$\omega$                  & 0.00281           & 0.00006 &   0.00269   & 0.00006	   \\
                          & ({\it  0.000014}) & ({\it 0.00002})  & ({\it 0.00020}) &  ({\it 0.00002}) \\
$a_1$                     & 0.99312 & 0.13495  & 0.00443  & 0.11422	    \\
                          & ({\it   0.02886}) & ({\it   0.01838}) & ({\it 0.00001}) & ({\it 0.01906}) \\
$b_1$                     & 0.00527   &  0.85318  & 0.02187  & 0.86526    \\
                          & ({\it 0.02182}) & ({\it 0.01752}) & ({\it 0.00001}) & ({\it 0.01704}) \\
$\sigma_{w}^2$            & 0.00043        &  0.00001  &  0.00007  &  0.00001\\
                          & ({\it 0.00041}) & ({\it 0.00001}) & ({\it 0.00016}) & ({\it 0.00003})  \\
$\beta_1^{(0)}$           & 0.22970        &  0.08714  & 0.15444	& -0.50465   \\
 ({\it initial})          & ({\it 1.32472}) & ({\it 0.39359}) &  ({\it 1.02795}) & ({\it 1.11751}) \\
\hline
\hline
$\max {\cal L}$ &  2688.9  & 2865.6 &  2710.5  & 2878.8\\
$\| \hat \rho_1(t) - \hat \beta_1(t)\|_{\infty}$ & 0.3045    &    0.3712 & 0.2988 & 0.3295  \\
\bottomrule
\end{tabular}
\end{table*}

The estimates calculated by the novel EKF-based approach for both the S\&P500 and DJIA series are summarized in Table~\ref{Table:estimate:AR1RW:GARCH}. Additionally, we perform a comparative study of the new method with the previously published KF-based algorithms utilized for evolving market efficiency tests in~\cite{2008:Posta,1999:Zalewska}. For a fair comparison, each adaptive filtering
strategy utilizes precisely the same sampled data, i.e. the S\&P500 and the DJIA  monthly returns, the same
initial system parameters and the same optimization method implemented. Having analyzed the results collected in Table~\ref{Table:estimate:AR1RW:GARCH}, we conclude that the estimation errors of the novel EKF-based method are smaller than for the classical KF procedure previously used in econometric literature for both return series under examination. As anticipated, the difference in their estimation quality is more significant than that obtained while estimating the time-varying AR(1) models with the homogeneous conditional variance assumption in Section~\ref{Section:AR}. This means that the EKF-based methodology provides a more accurate tracking algorithm for the level of market efficiency.

Finally, the estimated paths $\hat \beta_1(t_k)$ recovered by the novel EKF-based approach are illustrated by Figs.~\ref{beta1:General:graph} and~\ref{beta1:General:graph:DJIA} for the S\&P500 and the DJIA  monthly returns, respectively. For each sample examined, it is clearly seen that the trajectory $\left\{ \hat \beta_1(t_k)\right\}$ recovered by the filtering method adequately follows the sample autocorrelation coefficient $\left\{\hat \rho_1 (t_k)\right\}$ computed by the moving window method in Section~\ref{sec:data}, which is a proxy for an evolving market efficiency level.
 Further improvements in detecting the periods of anomalies might be achieved by designing the higher order estimation methods than the EKF technique.
 %For the models with sophisticated multiplicative noise term in the dynamic equation, the problem of deriving optimal nonlinear filters is of special interest because of a more accurate tracking method to be derived.

\section{Concluding remarks} \label{Section:conclusion}

In this paper, we proposed the nonlinear Bayesian estimation techniques for tracking a dynamics of weak form market efficiency. Taking into account a nonlinearity of the underlying models, the novel adaptive filtering framework is able to recover the time-varying AR coefficients from return series that correspond well to the sample autocorrelations computed by the moving window method. The interesting topics for a future research include the following issues: i) a further delicate modeling of the regression coefficients' dynamics as well as the conditional variance changes over time, and ii) improving a quality of tracking methods by designing the higher order estimators than the EKF-based methodology. A special attention should be paid to the derivative-free filtering techniques, which simplify their application in practice.

\section*{Acknowledgements}
The authors acknowledge the financial support of the Portuguese FCT~--- \emph{Funda\c{c}\~ao para a Ci\^encia e a Tecnologia}, through the projects UIDB/04621/2020 and UIDP/04621/2020 of CEMAT/IST-ID, Center for Computational and Stochastic Mathematics, Instituto Superior T\'ecnico, University of Lisbon 
and through the \emph{Scientific Employment Stimulus - 4th Edition} (CEEC-IND-4th edition) programme, grant number 2021.01450.CEECIND. They are also grateful to the anonymous referees for their valuable remarks and comments on the
paper.

%\bibliographystyle{elsarticle-num}
%\bibliography{BibTex_Library/books,%
%              BibTex_Library/KF_finance,%
%              BibTex_Library/KF_MLE,%
%              BibTex_Library/KF_Riccati,%
%              BibTex_Library/KFDiff_Riccati,%
%              BibTex_Library/Lit_finance}

\appendix
\section{Extended Kalman filtering and parameters estimation procedure}

The estimation problem considered in this paper is the unknown dynamic state and system parameters estimation of nonlinear discrete-time system given as follows~\cite[Section~5.1]{1970:Jazwinski:book}:
 \begin{align}
   x_{k+1} = & \phi(x_k,t_{k+1},t_k) + \Gamma(x_k,t_k)u_{k+1}, \quad k =0,1, \ldots, \label{eq:state:space1}
 \end{align}
where $x_k:= x(t_k)$ is an unknown $n$-vector to be estimated, $\phi(\cdot)$ is an $n$-vector function,  $\Gamma(\cdot)$ is $n \times q$, and $\{u_k, k=1,\ldots \}$ is an $q$-vector, white Gaussian sequence, $u_k \sim {\cal N}(0,Q_k)$. The distribution of the initial condition $x_0$ is assumed given, say $x_0 \sim {\cal N}(\bar x_0, P_0)$, and $x_0$ is independent of $\{u_k\}$. The noisy $m$-vector observations (measurements) $z_k$ to be  given by
 \begin{align}
   z_{k} = & h(x_k,t_k) + v_{k}, \quad k = 1,2 \ldots, \label{eq:state:space2}
 \end{align}
where $h(\cdot)$ is an $m$-vector function and $\{v_k: k=1,\ldots \}$ is an $m$-vector, Gaussian white sequence, $v_k \sim {\cal N}(0,R_k)$, $R_k>0$. For simplicity, $\{ u_k\}$ and $\{ v_k\}$ are assumed to be independent, and $\{ v_k\}$ is independent of $x_0$.

The estimation of GARCH models yields a more sophisticated case of non-additive noise in dynamic equation~\eqref{eq:state:space1}, i.e.
 \begin{align}
   x_{k+1} = & f(x_k,t_{k+1},t_k,u_{k+1}), \quad k =0,1, \ldots, \label{eq:state:space3}
 \end{align}
where $f(\cdot)$ is an $n$-vector function.

Clearly, the classical KF is not applicable for estimating such models and a proper solution has been found via the local linearization of nonlinear process and measurement equations around the filtering state estimate at each sampling time instant. This gave rise to a nonlinear state estimation method termed the {\it Extended Kalman Filter} (EKF). For decades, it was viewed as a simplest but successful and commonly used technique.

To derive the EKF equations, the discrete problem in~\eqref{eq:state:space3} is to be further linearized. For that, one expands its right-hand side function
$f(x_k,t_{k+1},t_k,u_{k+1})$ at the filtering state expectation
$\hat x_{k|k}$ of the state vector $x_{k}$ and around $\E{u_{k+1}}=0$. Thus, the related linear stochastic state-space system is derived~\cite{2016:SIAM:Kulikov}:
\begin{equation}\label{eq2.2}
x_{k}\approx\tilde x_{k}:=f\bigl(\hat x_{k|k},t_{k+1},t_k,0\bigr)+F_{k}(x_{k}-\hat x_{k|k})+G_{k}w_{k+1}
\end{equation}
where the Jacobian matrices are defined as follows: $F_{k} = \left.\partial f\bigl(x_{k},t_{k+1},t_k,0\bigr)/\partial x_k\right|_{\hat x_{k|k}}$,   $G_{k}=\left.\partial f\bigl(\hat x_{k|k},t_{k+1},t_k,u_{k+1}\bigr)/\partial u_{k+1}\right|_{0}$.

Next, the EKF equations are derived as usual, i.e.
\begin{align}
\hat x_{k+1|k}& = \E{\tilde x_k}=f\bigl(\hat x_{k|k},t_{k},t_{k+1},0\bigr), \label{state_new}\\
P_{k+1|k}&= \E{(\tilde x_k-\hat x_{k+1|k})(\tilde x_k-\hat x_{k+1|k})^\top} \nonumber \\
& = F_{k}P_{k|k}F_{k}^\top + G_{k}Q_{k}G_{k}^\top. \label{cov_new}
\end{align}

The measurement update of the EKF is summarized as follows~\cite[Theorem~8.1]{1970:Jazwinski:book}:
\begin{eqnarray}
K_k          & = & P_{k+1|k}H_k^{\top} R_{e,k}^{-1}, \; R_{e,k} = H_kP_{k+1|k}H_k^{\top} + R_k, \label{kf:f:K} \\
\hat x_{k+1|k+1} & = & \hat x_{k+1|k} + K_k e_k, \; e_k = z_k - h(\hat x_{k+1|k},t_k), \label{kf:f:X} \\
P_{k+1|k+1}      & = & (I -K_kH_k)P_{k+1|k} \label{kf:f:P}
\end{eqnarray}
where the matrix $H_k$ is defined by $H_k = \left.\partial h(x_k,t_k)/\partial x_k\right|_{\hat x_{k+1|k}}$.

\subsection{Filter initialization problem}

Since the state vector is a hidden stochastic process, information about the uncertainty of $x_0$ cannot be measured experimentally. Therefore, one needs to estimate the degree of uncertainty at the initial step. In general, if $\bar x_0$ is not close to $x_0$, then the filter's convergence to the correct estimate may be slow. It is important that if $P_0$ is chosen too small while $\bar x_0$ and $x_0$ differ significantly, then the KF as well as the extended KF may fail as reported, for instance, in chemical engineering literature~\cite{2005:Haseltine} and many other studies. The source of the divergence problem is in the incorrect initial filtering values when a small covariance implies a high confidence into the misleading initial state estimate when $\bar x_0$ is far from $x_0$. In this case, the filter might learn the wrong state too well and diverge~\cite{1970:Jazwinski:book}. It is commonly accepted in engineering and econometric literature that the filter should be initialized by $\Pi_0 = \delta I_{n}$ where $\delta  = \infty$ in case of no {\it a priori} information about the state available; e.g., see the discussion in~\cite{1990:Harvey,2012:Durbin:book,1987:Harvey:book,2015:Grewal:book}. However, the condition required yields a difficulty in practical implementation of  traditional {\it covariance}-type filtering algorithms because of the required setting $P_0=\infty$. As a result, the so-called {\it information}-type filters should be derived and implemented~\cite[p.~356-357]{2015:Grewal:book}. The information-form KF algorithms
process the information matrix, $\Lambda_{k|k}=P_{k|k}^{-1}$ instead, and hence the initial step is simple $\Lambda_{0} = P_0^{-1} = 0$. Some other advantages of the information-type KF implementation methods are discussed in~\cite[Section~7.7.3]{2015:Grewal:book}.

In contrast to the straightforward approach for resolving the filter initialization problem discussed above, some other powerful alternatives can be found in~\cite{1997:Koopman,2013:Schneider}. For the problem examined in this paper, we know upper and lower bounds on the initial state entries $\beta_{i,0}$, $i=1,\ldots,n$. This can be used for approximating the corresponding entries of $\bar x_0$ as suggested in~\cite{2013:Schneider}. More precisely, the hidden processes $\beta_{i,t}$, $i=1,\ldots,n$ are a proxy for the autocorrelation coefficients $\rho_{i,t}$ at lags $i=1,\ldots,n$. Hence, we conclude that $\beta_u = 1$ and $\beta_l = -1$. Next, the related entries of the initial state $\bar x_0$ are calculated as suggested in~\cite{2013:Schneider}:
\begin{align*}
\bar \beta_0 & = 0.5 (\beta_u + \beta_l), &  \bar \beta_0 & = [\beta_{1,0}, \ldots, \beta_{n,0}] = 0_{n},
\end{align*}
In summary, if the entries of the state vector $x_k = [x_1(t_k), x_2(t_k)]^{\top}:= [\beta_1(t_k), y(t_k)]^{\top}$, then the initial state is set to $x_0 = [\bar \beta_0, y(t_1)]^{\top}$ where $\bar \beta_0 = 0$.

Next, the corresponding diagonal values of the filter error covariance matrix $P_0$ might be defined as follows~\cite{2013:Schneider}:
\begin{align*}
\alpha & = \bar \beta_0-\beta_0 = \max\{|\beta_u - \bar \beta_0|, |\beta_l - \bar \beta_0|\}, &  P_0 & = \alpha^2 \cdot I_n = I_n
\end{align*}
that is $P_0=I_n$ in case of zero initial value $\bar \beta_0$ utilized.

Finally, the most general way is to consider the initial values $\beta_{i,0}$, $i=1,\ldots,n$ as the extra unknown system parameters and next estimate them in a systematic manner discussed below.

\subsection{System parameters estimation and filter's tuning}

The process noise covariance $Q$ and the measurement noise covariance matrix $R$ are often assumed to be constant over time, but they are usually unknown. %Consequently, one needs to select these matrices prior to filtering. The measurement noise covariance $R$ might be estimated from the actual measurements by taking some off-line sample measurements in order to determine the variance of the measurement noise. The determination of the process noise covariance $Q$ is generally more difficult since we usually do not have the ability to observe the process to be estimated.
It is preferable to define the unknown system parameters, say $\theta$, in a systematic manner rather than {\it ad hoc} by trial and error approach or by setting arbitrarily chosen values. This means that the state-space model examined is parameterized. The most frequently used estimation approach is the method of maximum likelihood, i.e. the model is fitted to the data by maximizing the log likelihood function given as follows~\cite{1989:Harvey:book,1965:Schweppe} (without a constant term):
\begin{align}
{\cal L}(\theta) = -\frac{1}{2}\sum \limits_{k=1}^N \ln \det R_{e,k}\! - \!\frac{1}{2}\sum \limits_{k=1}^N \left\{e_k^\top R_{e,k}^{-1}e_k \right\} \label{eq:llf}
\end{align}
where the residual $e_k$ and covariance $R_{e,k}$ come from the filtering algorithm, i.e. they are computed from formulas~\eqref{kf:f:X} and~\eqref{kf:f:K}, respectively.

In summary, any adaptive filtering scheme includes the filtering method for computing the performance index in~\eqref{eq:llf} and  optimization method for finding the optimal parameters value $\hat \theta_{MLE}^*$. In this paper, we use the built-in MATLAB function \verb"fmincon" for optimization purposes with the tolerance $e_{\theta} = 10^{-6}$ on the parameters value and with $e_{log LF} = 10^{-9}$ on the function in the stopping criterion applied. It is also worth noting here that  econometric models are often estimated by the Expectation-Maximization (EM) algorithm; e.g. see~\cite{1982:Shumway,1993:Koopman,2007:Elliott}, and many others. In practice, the EM algorithm might be combined with
gradient-based methods for a faster convergence as discussed in~\cite{2020:Galka}.

\end{document}